\documentclass[10pt]{article}
\usepackage{tikz}
\tikzset{dashdot/.style={dash pattern=on 4pt off 1pt on 1pt off 1pt}}
\usepackage{multirow}
\usepackage{url}
\usepackage{amsthm,amssymb}
\usepackage{amsfonts}
\usepackage{verbatim}
\usepackage{setspace}
\usepackage{float}
\theoremstyle{plain}
\newtheorem{theorem}{Theorem}
\newtheorem{lemma}[theorem]{Lemma}
\newtheorem{proposition}[theorem]{Proposition}
\newtheorem{corollary}[theorem]{Corollary}

\theoremstyle{definition}

\DeclareSymbolFont{AMSb}{U}{msb}{m}{n}
\DeclareMathSymbol{\N}{\mathbin}{AMSb}{"4E}
\DeclareMathSymbol{\Z}{\mathbin}{AMSb}{"5A}
\DeclareMathSymbol{\R}{\mathbin}{AMSb}{"52}
\DeclareMathSymbol{\Q}{\mathbin}{AMSb}{"51}
\DeclareMathSymbol{\I}{\mathbin}{AMSb}{"49}
\DeclareMathSymbol{\C}{\mathbin}{AMSb}{"43}

\title{Lehmer's conjecture for matrices over the ring of integers of some imaginary quadratic fields.}
\author{Graeme Taylor\footnote{Heilbronn Institute for Mathematical Research, University of Bristol, School of Mathematics, Howard House, Queens Avenue, Bristol BS8 1SN, UK}}
\date{}
\begin{document}
\maketitle
\begin{abstract} \noindent Let $R=\mathcal{O}_{\Q(\sqrt{d})}$ for $d<0$, squarefree, $d\neq -1,-3$. We prove Lehmer's conjecture for associated reciprocal polynomials of $R$-matrices; that is, any noncyclotomic $R$-matrix has Mahler measure at least $\lambda_0=1.176\ldots$. \end{abstract}

\section{Introduction}

Given a monic polynomial $P(z)=\displaystyle\prod_{i=1}^{d}(z-\alpha_i)\in\Z[z]$, the \emph{Mahler Measure} $M(P)$ is given by \[M(P):=\prod_{i=1}^d \max{(1,|\alpha_i|)}=\prod_{|\alpha_i|>1} |\alpha_i|\] Clearly $M(P)\ge1$; by a result of Kronecker (\cite{Kron}) $M(P)=1$ if and only if $\pm P$ is the product of a cyclotomic polynomial\footnote{Following Boyd \cite{Boyd}, we will use `cyclotomic' to refer to any polynomial for which all roots are roots of unity, rather than just the irreducible examples $\Phi_n$.} and a power of $z$. For a monic integer polynomial with $M(P)>1$, Lehmer asked (in \cite{Lehmer}) whether $M(P)$ could be arbitrarily close to $1$. This is now known as \emph{Lehmer's Problem}; the negative result - that there is some $\lambda>1$ such that $M(P)>1 \Rightarrow M(P)\ge\lambda$ - is sometimes referred to as \emph{Lehmer's Conjecture}. In \cite{Lehmer}, Lehmer exhibited the polynomial \[z^{10}+z^9-z^7-z^6-z^5-z^4-z^3+z+1\] with Mahler measure $\lambda_0=1.176\ldots$; no noncyclotomic monic integer polynomial with lesser Mahler measure has been found since.

For a monic polynomial $g\in\Z[x]$ of degree $n$, we define its \emph{associated reciprocal polynomial} to be $z^ng(z+1/z)$ which is a monic reciprocal polynomial of degree $2n$. For $R=\mathcal{O}_{\Q(\sqrt{d})}$ with $d<0$, let $A$ be an $n\times n$ Hermitian $R$-matrix and denote by $R_A(z)$ the associated reciprocal polynomial of its characteristic polynomial $\chi_A(x)=\mbox{det}(xI-A)$. Further, define $M(A)$, \emph{the Mahler measure of $A$}, to be $M(R_A(z))$. Since $\chi_A\in\Z[x]$ with all roots real, if $A$ has spectral radius at most $2$ then $R_A(z)$ is cyclotomic; we describe such an $A$ as a \emph{cyclotomic matrix}. A classification for cyclotomic integer symmetric matrices was given by McKee and Smyth in \cite{McSm1}; in \cite{McSm2} they subsequently proved that any noncyclotomic integer symmetric matrix has Mahler measure at least $\lambda_0$. However, they were also able to demonstrate the existence of noncyclotomic reciprocal polynomials $f\in\Z[z]$ such that no integer symmetric matrix $A$ satisfies $M(A)=M(f)$. These `missing' Mahler measures motivated the extension to broader classes of matrices: the cyclotomic $R$-matrices for squarefree $d<0$ are found in \cite{Ta2}  ($d\neq -1,-3$) and \cite{Gr1} ($d=-1,-3$). In this work we build upon the results of \cite{McSm2} and \cite{Ta2} to prove:

\begin{theorem}\label{lehmer271115} Let $A$ be a Hermitian $\mathcal{O}_{\Q(\sqrt{d})}$-matrix for squarefree $d<0$, $d\neq -1,-3$. Then \[M(A)=1 \mbox{ or } M(A)\ge\lambda_0\]
\end{theorem}

In Section 2, we survey the results for cyclotomic matrices and introduce corresponding graph structures. We will then demonstrate that any noncyclotomic $R$-matrix with `small' Mahler measure is in fact an integer symmetric matrix. To do so, we eliminate the possibility of large norm entries (Section 3); search for examples with at most ten vertices (Section 4); and prove that this search was exhaustive by showing there can be no larger examples (Section 5). These results are assembled into proofs of Theorem \ref{lehmer271115} for each $d$ in Section 6.

\section{Cyclotomic Matrices and Graphs}
\emph{Throughout, we assume $R=\mathcal{O}_{\Q(\sqrt{d})}$ for $d<0$ squarefree and $d\neq -1,-3$, and that (unless otherwise stated) all matrices are Hermitian and all graphs are connected.}

\subsection{Cyclotomic Matrices}
If $A$ is a block diagonal matrix, then its list of eigenvalues is the union of the lists of the eigenvalues of the blocks. If there is a reordering of the rows (and columns) of $A$ such that it has block diagonal form with more than one block, then $A$ will be called \emph{decomposable}; if there is no such reordering, $A$ is called \emph{indecomposable}. Clearly any decomposable cyclotomic matrix decomposes into cyclotomic blocks. But a much stronger result holds:

\begin{theorem}[Cauchy Interlacing Theorem \cite{Cauchy}]\label{cauchyinterlacetheorem}
Let $A$ be a Hermitian $n\times n$ matrix with eigenvalues $\lambda_1\le\lambda_2\le\cdots\le\lambda_n$. \\Let $B$ be obtained from $A$ by deleting row $i$ and column $i$ from $A$. \\Then the eigenvalues $\mu_1\le\cdots\le\mu_{n-1}$ of $B$ interlace with those of $A$: that is, \[ \lambda_1\le\mu_1\le\lambda_2\le\mu_2\le\cdots\le\lambda_{n-1}\le\mu_{n-1}\le\lambda_n .\]
\end{theorem}

Thus if $A$ is cyclotomic, so is any $B$ obtained by successively deleting a series of rows and corresponding columns from $A$. We describe such a $B$ as being \emph{contained} in $A$. If an indecomposable cyclotomic matrix $A$ is not contained in a strictly larger indecomposable cyclotomic matrix, then we call $A$ \emph{maximal}.

Additionally, an equivalence relation on cyclotomic $R$-matrices can be defined as follows. Let $O_n(\Z)$ denote the orthogonal group of $n\times n$ signed permutation matrices. Conjugation of a cyclotomic matrix by a matrix from this group gives another matrix with the same eigenvalues, which is thus also cyclotomic. Cyclotomic matrices $A,A'$ related in this way are described as \emph{strongly equivalent}; indecomposable cyclotomic matrices $A$ and $A'$ are then considered equivalent if $A'$ is strongly equivalent to any of $A$, $-A$, $\overline{A}$ or $-\overline{A}$.

We note the following constraint on entries of cyclotomic $R$-matrices:

\begin{lemma}[\cite{Ta2} Lemma 5]\label{norm4max} Let $A$ be a cyclotomic $R$-matrix. Then for any entry $A_{ij}$ of $A$, \[ | A_{ij}A_{ji}| \le 4.\]\end{lemma}

I.e., if $i=j$ then $|A_{ii}|\le2$, and if $i\neq j$ then $\mbox{Norm}(A_{ij})\le4$.  In fact (see \cite{Ta2}), we have that for any indecomposable cyclotomic $R$-matrix $A$, $A_{ii}\in\{0,1,-1\}$ unless $A$ is the $1\times1$ matrix $(2)$ or $(-2)$.

For a given $R$ and $n\ge1$, define $\mathcal{L}_n=\{x\in R\,|\,x\overline{x}=n\}$.  Then if $A$ is a cyclotomic Hermitian matrix with all entries from $R$, by Lemma \ref{norm4max} $A$ is an $\mathcal{L}$-matrix for \[ \mathcal{L}:= \{0\} \cup \mathcal{L}_1 \cup \mathcal{L}_2 \cup \mathcal{L}_3 \cup \mathcal{L}_4.\] For convenience, we also define \[ \mathcal{L}':= \{0\} \cup \mathcal{L}_1 \cup \mathcal{L}_2.\]

\begin{corollary} For a squarefree, negative $d\not\in\{-1,-2,-3,-7,-11,-15\}$, $\mathcal{L}=\{0,\pm1,\pm2\}$ and thus any cyclotomic Hermitian $\mathcal{L}$-matrix is an integer symmetric matrix. \end{corollary}

\subsection{Cyclotomic $\mathcal{L}$-graphs}

We may construct an $n$-vertex $\mathcal{L}$-graph $G$ from an $n\times n$ Hermitian $\mathcal{L}$-matrix by specifying nonzero entries of $A$ as edge- or vertex weights for $G$. 

For each vertex $i$ of $G$ corresponding to the diagonal entry $A_{ii}$ of $A$, $i$ can either be \emph{neutral} ($A_{ii}=0$), \emph{positive} ($A_{ii}=1$) or \emph{negative} ($A_{ii}=-1$); we indicate these visually as $\begin{tikzpicture}
\node (a) at (0,0) [fill=black,draw,shape=circle] {};
\end{tikzpicture}$,
$\begin{tikzpicture}
\node (a) at (0,0) [fill=white,inner sep=0pt,draw,shape=circle] {$+$};
\end{tikzpicture}$ and
$\begin{tikzpicture}
\node (a) at (0,0) [fill=white,inner sep=0pt,draw,shape=circle] {$-$};
\end{tikzpicture}$ respectively, and describe the latter two as \emph{charged} vertices. 

For each $i<j$ such that $A_{ij}=x \in\mathcal{L}\backslash\{0\}$ we introduce an edge with label $x$ from vertex $i$ to vertex $j$. We define the \emph{weight} of an edge to be the norm of $x$, so a \emph{weight $n$ edge} is one with a label from $\mathcal{L}_n$. For each of our chosen $R$, the only possible weight one edge labels are $+1$ and $-1$, so we may speak of the \emph{sign} of such an edge: edges with a positive sign will be drawn as $\begin{tikzpicture}
\node (a) at (0,0) {};
\node (b) at (1,0) {};
\draw (a) -- (b);
\end{tikzpicture}$ and edges with a negative sign will be drawn as $\begin{tikzpicture}
\node (a) at (0,0) {};
\node (b) at (1,0) {};
\draw [dotted,thick] (a) -- (b);
\end{tikzpicture}$. For higher weight edges, we need specify the label, but can indicate the weight visually by denoting an edge label $\omega$ from $\mathcal{L}_2$, $\mathcal{L}_3$ or $\mathcal{L}_4$ by $\begin{tikzpicture}[scale=1]
\node (a) at (0,0) [fill=black,draw,shape=circle] {};
\node (b) at (1,0) [fill=black,draw,shape=circle] {};
\draw [double] (a) -- node[above] {\scriptsize $\omega$} (b);\end{tikzpicture}$, $\begin{tikzpicture}[scale=1]
\node (a) at (0,0) [fill=black,draw,shape=circle] {};
\node (b) at (1,0) [fill=black,draw,shape=circle] {};
\draw (a.0) -- (b.180) (a.340) -- (b.200) (a.20) -- node[above] {\scriptsize$\omega$} (b.160);\end{tikzpicture}$ or $\begin{tikzpicture}[scale=1]
\node (a) at (0,0) [fill=black,draw,shape=circle] {};
\node (b) at (1,0) [fill=black,draw,shape=circle] {};
\draw [double] (a.15) -- node[above] {\scriptsize $\omega$} (b.165);
\draw [double] (a.345) -- (b.195);\end{tikzpicture}$ respectively. 

If $A$ is in fact a symmetric $\{0,1,-1\}$ matrix with only zeros on the diagonal, then its $\mathcal{L}$-graph is a \emph{signed graph} as in \cite{CST}, \cite{Zas}: allowing only positive entries gives the usual identification of adjacency matrices and graphs; whilst allowing nonzero entries on the diagonal gives the \emph{charged signed graphs} of \cite{McSm1}, \cite{McSm2}.

An $\mathcal{L}$-graph $G$ is described as cyclotomic if its corresponding $\mathcal{L}$-matrix $A$ is cyclotomic; the Mahler measure of $G$ is that of $A$ (i.e., of $R_A(z)$), and $\mathcal{L}$-graphs $G,G'$ are (strongly) equivalent if and only if their $\mathcal{L}$-matrices $A,A'$ are. An $\mathcal{L}$-graph $G$ is connected if and only if its corresponding $\mathcal{L}$-matrix is indecomposable. If a cyclotomic $\mathcal{L}$-matrix $A'$ is contained in $A$ then its corresponding $\mathcal{L}$-graph $G'$ is an induced subgraph of $G$ corresponding to $A$; thus a maximal cyclotomic $\mathcal{L}$-graph is connected yet not an induced subgraph of any strictly larger connected cyclotomic $\mathcal{L}$-graph.

The equivalence relation on $\mathcal{L}$-matrices has the following interpretation for $\mathcal{L}$-graphs. $O_n(\Z)$ is generated by matrices of the form $\mbox{diag}(1,1,\ldots,1,-1,1,\ldots,1)$ and permutation matrices. Conjugation by the former has the effect of negating the signs of all edges incident at some vertex $v$; following \cite{CST} this is described as \emph{switching at $v$}. Conjugation by a permutation matrix merely permutes vertex labels; for charged signed graphs, we may therefore omit vertex labels.
 
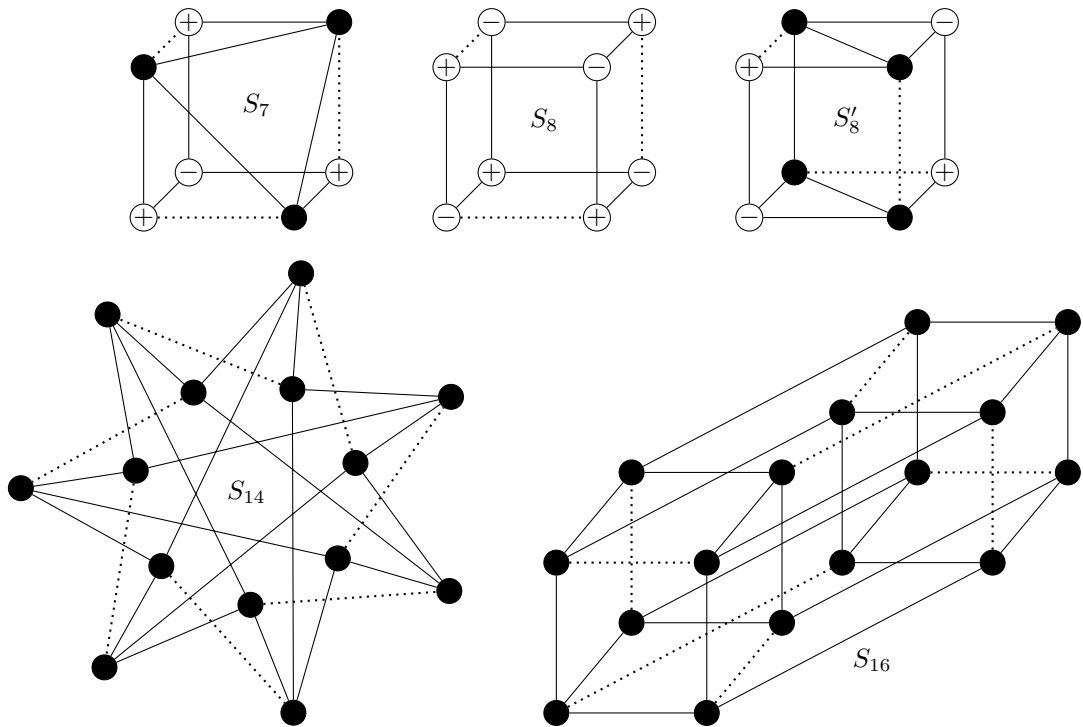
\begin{figure}[H]

\[
\mbox{$\begin{tikzpicture}[scale=2]
\node (1) at (0,0) [fill=white,inner sep=0pt,draw,shape=circle] {$+$};
\node (2) at (1,0) [fill=black,draw,shape=circle] {};
\node (3) at (0,1) [fill=black,draw,shape=circle] {};
\draw (1) -- (3);
\draw [dotted,thick] (1) -- (2);
\node (1b) at (0.3,0.3) [fill=white,inner sep=0pt,draw,shape=circle] {$-$};
\node (2b) at (1.3,0.3) [fill=white,inner sep=0pt,draw,shape=circle] {$+$};
\node (3b) at (0.3,1.3) [fill=white,inner sep=0pt,draw,shape=circle] {$+$};
\node (4b) at (1.3,1.3) [fill=black,draw,shape=circle] {};
\draw (2b) -- (1b) -- (3b) -- (4b);
\draw [dotted,thick] (2b) -- (4b);
\draw (1) -- (1b) (2) -- (2b);
\draw [dotted,thick] (3) -- (3b);
\draw (2) -- (3) -- (4b) -- (2);
\node (name) at (0.75,0.75) {$S_7$};
\end{tikzpicture}$}
\hspace{3em}
\mbox{$\begin{tikzpicture}[scale=2]
\node (1) at (0,0) [fill=white,inner sep=0pt,draw,shape=circle] {$-$};
\node (2) at (1,0) [fill=white,inner sep=0pt,draw,shape=circle] {$+$};
\node (3) at (0,1) [fill=white,inner sep=0pt,draw,shape=circle] {$+$};
\node (4) at (1,1) [fill=white,inner sep=0pt,draw,shape=circle] {$-$};
\draw (1) -- (3) -- (4) -- (2);
\draw [dotted,thick] (1) -- (2);
\node (1b) at (0.3,0.3) [fill=white,inner sep=0pt,draw,shape=circle] {$+$};
\node (2b) at (1.3,0.3) [fill=white,inner sep=0pt,draw,shape=circle] {$-$};
\node (3b) at (0.3,1.3) [fill=white,inner sep=0pt,draw,shape=circle] {$-$};
\node (4b) at (1.3,1.3) [fill=white,inner sep=0pt,draw,shape=circle] {$+$};
\draw (2b) -- (1b) -- (3b) -- (4b);
\draw [dotted,thick] (2b) -- (4b);
\draw (1) -- (1b) (2) -- (2b) (4) -- (4b);
\draw [dotted,thick] (3) -- (3b);
\node (name) at (0.65,0.65) {$S_8$};
\end{tikzpicture}$}
\hspace{3em}
\mbox{$\begin{tikzpicture}[scale=2]
\node (1) at (0,0) [fill=white,inner sep=0pt,draw,shape=circle] {$-$};
\node (2) at (1,0) [fill=black,draw,shape=circle] {};
\node (3) at (0,1) [fill=white,inner sep=0pt,draw,shape=circle] {$+$};
\node (4) at (1,1) [fill=black,draw,shape=circle] {};
\draw (2) -- (1) -- (3) -- (4);
\draw [dotted,thick] (4) -- (2);
\node (1b) at (0.3,0.3) [fill=black,draw,shape=circle] {};
\node (2b) at (1.3,0.3) [fill=white,inner sep=0pt,draw,shape=circle] {$+$};
\node (3b) at (0.3,1.3) [fill=black,draw,shape=circle] {};
\node (4b) at (1.3,1.3) [fill=white,inner sep=0pt,draw,shape=circle] {$-$};
\draw (2b) -- (4b) -- (3b) -- (1b);
\draw [dotted,thick] (1b) -- (2b);
\draw (1) -- (1b) (2) -- (2b) (4) -- (4b) (2) -- (1b) (4) -- (3b);
\draw [dotted,thick] (3) -- (3b);
\node (name) at (0.65,0.65) {$S_8'$};
\end{tikzpicture}$}
\]
	\[
\mbox{$\begin{tikzpicture}[scale=1.5]
\node (1) at (15:1) [fill=black,draw,shape=circle] {};
\node (2) at (66:1) [fill=black,draw,shape=circle] {};
\node (3) at (118:1) [fill=black,draw,shape=circle] {};
\node (4) at (169:1) [fill=black,draw,shape=circle] {};
\node (5) at (221:1) [fill=black,draw,shape=circle] {};
\node (6) at (272:1) [fill=black,draw,shape=circle] {};
\node (7) at (324:1) [fill=black,draw,shape=circle] {};
\node (1a) at (25:2) [fill=black,draw,shape=circle] {};
\node (2a) at (76:2) [fill=black,draw,shape=circle] {};
\node (3a) at (128:2) [fill=black,draw,shape=circle] {};
\node (4a) at (179:2) [fill=black,draw,shape=circle] {};
\node (5a) at (231:2) [fill=black,draw,shape=circle] {};
\node (6a) at (282:2) [fill=black,draw,shape=circle] {};
\node (7a) at (334:2) [fill=black,draw,shape=circle] {};
\draw (1) -- (1a) -- (2) -- (2a) -- (3) -- (3a) -- (4) -- (4a) -- (5) -- (5a) -- (6) -- (6a) -- (7) -- (7a) -- (1);
\draw (1) -- (5a) (2) -- (6a) (3) -- (7a) (4) -- (1a) (5) -- (2a) (6) -- (3a) (7) -- (4a);
\draw [dotted,thick] (1) -- (2a) (2) -- (3a) (3) -- (4a) (4) -- (5a) (5) -- (6a) (6) -- (7a) (7) -- (1a);
\node (label) at (0,0) {$S_{14}$};
\end{tikzpicture}$}
\hspace{3em}
\mbox{$\begin{tikzpicture}
\node (1) at (0,0) [fill=black,draw,shape=circle] {};
\node (2) at (2,0) [fill=black,draw,shape=circle] {};
\node (3) at (0,2) [fill=black,draw,shape=circle] {};
\node (4) at (2,2) [fill=black,draw,shape=circle] {};
\draw (3) -- (1) -- (2) -- (4);
\draw [dotted,thick] (3) -- (4);
\node (1b) at (1,1.2) [fill=black,draw,shape=circle] {};
\node (2b) at (3,1.2) [fill=black,draw,shape=circle] {};
\node (3b) at (1,3.2) [fill=black,draw,shape=circle] {};
\node (4b) at (3,3.2) [fill=black,draw,shape=circle] {};
\draw [dotted,thick] (1b) -- (3b);
\draw (1b) -- (2b) -- (4b) -- (3b);
\draw (1) -- (1b) (3) -- (3b) (4) -- (4b);
\draw [dotted,thick] (2) -- (2b);
\node (5) at (3.8,2) [fill=black,draw,shape=circle] {};
\node (6) at (5.8,2) [fill=black,draw,shape=circle] {};
\node (7) at (3.8,4) [fill=black,draw,shape=circle] {};
\node (8) at (5.8,4) [fill=black,draw,shape=circle] {};
\node (5b) at (4.8,3.2) [fill=black,draw,shape=circle] {};
\node (6b) at (6.8,3.2) [fill=black,draw,shape=circle] {};
\node (7b) at (4.8,5.2) [fill=black,draw,shape=circle] {};
\node (8b) at (6.8,5.2) [fill=black,draw,shape=circle] {};
\draw (6) -- (5) -- (7) -- (8);
\draw [dotted,thick](6) -- (8);
\draw (5b) -- (7b) -- (8b) -- (6b);
\draw [dotted,thick] (5b) -- (6b);
\draw (5) -- (5b) (6) -- (6b) (8) -- (8b);
\draw [dotted,thick] (7) -- (7b);
\draw (2) -- (6) (3) -- (7) (4) -- (8) (1b) -- (5b) (2b) -- (6b) (3b) -- (7b);
\draw [dotted,thick] (1) -- (5) (4b) -- (8b);
\node (label) at (4.2,0.7) {$S_{16}$};
\end{tikzpicture}$}
\]
\caption{The sporadic maximal cyclotomic charged signed graphs $S_7$, $S_8$, $S_8'$, $S_{14}$ and $S_{16}$.}
\label{sporadicCSGfigure}
\end{figure}

\begin{figure}[H]
\[\begin{tikzpicture}[scale=1.5]
\node (a) at (0,0) [fill=black,draw,shape=circle] {};
\node (b) at (1,0) [fill=black,draw,shape=circle] {};
\draw [double] (a.15) -- node[above] {\scriptsize $2$} (b.165);
\draw [double] (a.345) -- (b.195);
\node (a) at (2,0) [fill=black,draw,shape=circle] {};
\node (b) at (3,0) [fill=black,draw,shape=circle] {};
\draw [double] (a.15) -- node[above] {\scriptsize $\omega$} (b.165);
\draw [double] (a.345) -- (b.195);
\node (S2label) at (0.5,-0.5) {$S_2$};
\node (S2label1) at (2.5,-0.5) {$S_2^*$};
\end{tikzpicture}\]
\caption{The sporadic maximal cyclotomic $\mathcal{L}$-graphs $S_2$ and $S_2^*$ \emph{(Where $\omega=3/2+\sqrt{-7}/2, 1/2+\sqrt{-15}/2$ for $d=-7,-15$ respectively.)}}
\label{weight4figure}
\end{figure}
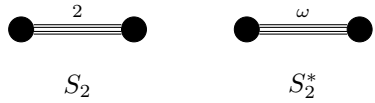

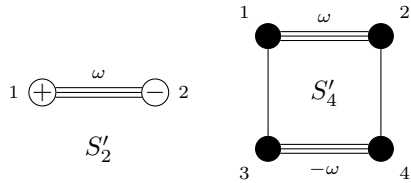
\begin{figure}[H]
\[\begin{tikzpicture}[scale=1.5]
\node (a) at (0,0) [fill=white,inner sep=0pt,draw,circle,label=left:{\scriptsize$1$}] {$+$};
\node (b) at (1,0) [fill=white,inner sep=0pt,draw,circle,label=right:{\scriptsize$2$}] {$-$};
\draw (a.0) -- (b.180) (a.340) -- (b.200) (a.20) -- node[above] {\scriptsize$\omega$} (b.160);
\node (S2label) at (0.5,-0.5) {$S_2'$};

\node (a) at (2,0.5) [fill=black,draw,circle,label=above left:{\scriptsize$1$}] {};
\node (b) at (3,0.5) [fill=black,draw,circle,label=above right:{\scriptsize$2$}] {};
\draw (a.0) -- (b.180) (a.340) -- (b.200) (a.20) -- node[above] {\scriptsize$\omega$} (b.160);
\node (c) at (2,-0.5) [fill=black,draw,circle,label=below left:{\scriptsize$3$}] {};
\node (d) at (3,-0.5) [fill=black,draw,circle,label=below right:{\scriptsize$4$}] {};
\draw (c.0) -- (d.180) (c.340) -- node[below] {\scriptsize$-\omega$} (d.200) (c.20) -- (d.160);
\draw (a) -- (c) (b) -- (d);
\node (S4label) at (2.5,0) {$S_4'$};
\end{tikzpicture}
\]
\caption{The sporadic maximal cyclotomic $\mathcal{L}$-graphs $S_2'$ and $S_4'$ \emph{(Where $\omega=1+\sqrt{-2}, 1/2+\sqrt{-11}/2$ for $d=-2,-11$ respectively.)}}
\label{weight3figure}
\end{figure}

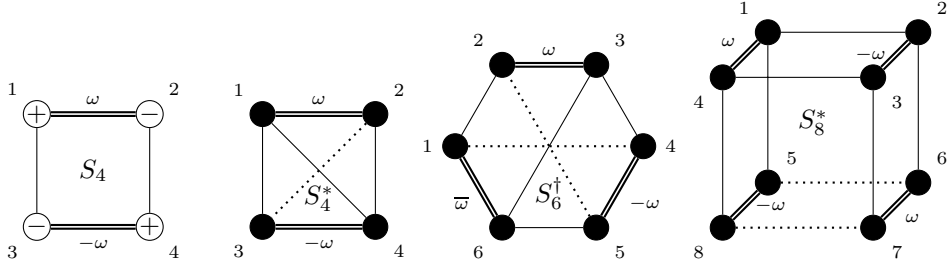
\begin{figure}[H]
\[\mbox{$\begin{tikzpicture}[scale=1.5]
\node (a) at (0,0) [fill=white,inner sep=0pt,draw,circle,label=above left:{\scriptsize$1$}] {$+$};
\node (b) at (1,0) [fill=white,inner sep=0pt,draw,circle,label=above right:{\scriptsize$2$}] {$-$};
\node (c) at (0,-1) [fill=white,inner sep=0pt,draw,circle,label=below left:{\scriptsize$3$}] {$-$};
\node (d) at (1,-1) [fill=white,inner sep=0pt,draw,circle,label=below right:{\scriptsize$4$}] {$+$};
\draw [double,thick] (a) -- node[above] {\scriptsize$\omega$} (b) (c) -- node[below] {\scriptsize$-\omega$} (d);
\draw (a) -- (c) (b) -- (d);
\node (S4label) at (0.5,-0.5) {$S_4$};
\node (a) at (2,0) [fill=black,draw,circle,label=above left:{\scriptsize$1$}] {};
\node (b) at (3,0) [fill=black,draw,circle,label=above right:{\scriptsize$2$}] {};
\node (c) at (2,-1) [fill=black,draw,circle,label=below left:{\scriptsize$3$}] {};
\node (d) at (3,-1) [fill=black,draw,circle,label=below right:{\scriptsize$4$}] {};
\draw [double,thick] (a) -- node[above] {\scriptsize$\omega$} (b);
\draw [double,thick] (c) -- node[below] {\scriptsize$-\omega$} (d);
\draw (a) -- (c) (b) -- (d) (a) -- (d);
\draw [dotted,thick] (b) -- (c);
\node (S4label) at (2.5,-0.75) {$S_4^*$};
\end{tikzpicture}$}
\mbox{$\begin{tikzpicture}[scale=1.25]
\node (a) at (0:1) [fill=black,draw,circle,label=right:{\scriptsize$4$}] {};
\node (b) at (60:1) [fill=black,draw,circle,label=above right:{\scriptsize$3$}] {};
\node (c) at (120:1) [fill=black,draw,circle,label=above left:{\scriptsize$2$}] {};
\node (d) at (180:1) [fill=black,draw,circle,label=left:{\scriptsize$1$}] {};
\node (e) at (240:1) [fill=black,draw,circle,label=below left:{\scriptsize$6$}] {};
\node (f) at (300:1) [fill=black,draw,circle,label=below right:{\scriptsize$5$}] {};
\draw (d) -- (c) (a) -- (b) (e) -- (f) (b) -- (e) ;
\draw [dotted,thick] (f) -- (c) (d) -- (a);
\draw [double,thick] (d) -- node[below left] {\scriptsize$\overline{\omega}$} (e) (f) -- node[below right] {\scriptsize$-\omega$} (a);
\draw [double,thick] (b) -- node[above] {\scriptsize$\omega$} (c);
\node (S6label) at (0,-0.5) {$S_6^\dag$};
\end{tikzpicture}$}
\mbox{$\begin{tikzpicture}[scale=2]
\node (1) at (2,0) [fill=black,draw,shape=circle,label=below left:{\scriptsize$8$}] {};
\node (2) at (3,0) [fill=black,draw,shape=circle,label=below right:{\scriptsize$7$}] {};
\node (3) at (2,1) [fill=black,draw,shape=circle,label=below left:{\scriptsize$4$}] {};
\node (4) at (3,1) [fill=black,draw,shape=circle,label=below right:{\scriptsize$3$}] {};
\node (1b) at (2.3,0.3) [fill=black,draw,shape=circle,label=above right:{\scriptsize$5$}] {};
\node (2b) at (3.3,0.3) [fill=black,draw,shape=circle,label=above right:{\scriptsize$6$}] {};
\node (3b) at (2.3,1.3) [fill=black,draw,shape=circle,label=above left:{\scriptsize$1$}] {};
\node (4b) at (3.3,1.3) [fill=black,draw,shape=circle,label=above right:{\scriptsize$2$}] {};
\draw [dotted,thick] (1) -- (2) (1b) -- (2b);
\draw (1) -- (3) (2) -- (4) (1b) -- (3b) (2b) -- (4b) (3) -- (4) (3b) -- (4b);
\draw [double,thick] (1) -- node[right] {\scriptsize$-\omega$} (1b) (2) -- node[right,yshift=-0.5em] {\scriptsize$\omega$} (2b) (4) -- node[left] {\scriptsize$-\omega$} (4b) (3) -- node[left,yshift=0.5em] {\scriptsize$\omega$} (3b);
\node (S8label) at (2.6,0.7) {$S_8^*$};
\end{tikzpicture}$}
\]

\caption{The sporadic maximal cyclotomic $\mathcal{L}'$-graphs $S_4$, $S_4^*$ ($d=-2$ only), $S_6^\dag$ ($d=-7$ only) and $S_8^*$ \emph{(Where $\omega=\sqrt{-2},1/2+\sqrt{-7}/2$ for $d=-2,-7$ respectively.)}}
\label{Sporadicw2figure}
\end{figure}

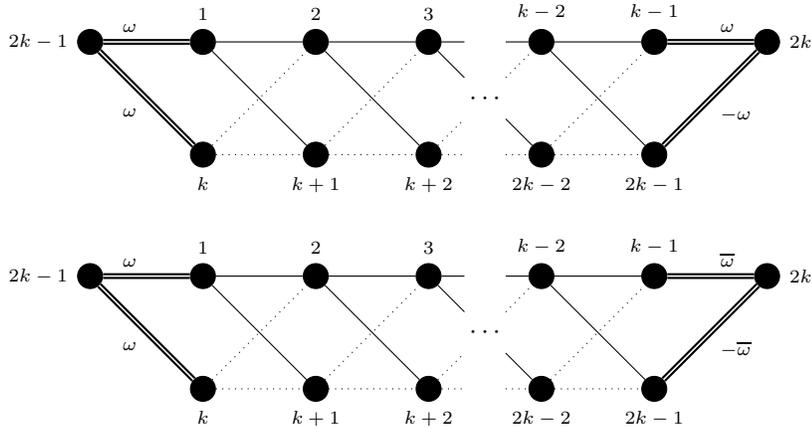
\begin{figure}[H]
\[
\begin{tikzpicture}[scale=1.5]
\node (1) at (0,0) [fill=black,draw,circle,label=above:{\scriptsize$1$}] {};
\node (2) at (1,0) [fill=black,draw,circle,label=above:{\scriptsize$2$}] {};
\node (3) at (2,0) [fill=black,draw,circle,label=above:{\scriptsize$3$}] {};
\node (3a) at (2.4,0) {};
\node (3b) at (2.4,-0.4) {};
\draw (1) -- (2) -- (3) -- (3a); 
\draw (3) -- (3b);
\node (L+1) at (0,-1) [fill=black,draw,circle,label=below:{\scriptsize$k$}] {};
\node (L+2) at (1,-1) [fill=black,draw,circle,label=below:{\scriptsize$k+1$}] {};
\node (L+3) at (2,-1) [fill=black,draw,circle,label=below:{\scriptsize$k+2$}] {};
\node (k+3a) at (2.4,-1) {};
\node (k+3b) at (2.4,-0.6) {};
\draw [dotted] (L+1) -- (L+2) -- (L+3) -- (k+3a);
\draw [dotted] (L+3) -- (k+3b);
\draw (1) -- (L+2) (2) -- (L+3);
\draw [dotted] (L+1) -- (2) (L+2) -- (3);
\node (spacer) at (2.5,-0.5) {$\cdots$};
\node (k-1b) at (2.6,-0.4) {};
\node (2k-1b) at (2.6,-0.6) {};
\node (k-1a) at (2.6,0) {};
\node (L-1) at (3,0) [fill=black,draw,circle,label=above:{\scriptsize$k-2$}] {};
\node (L) at (4,0) [fill=black,draw,circle,label=above:{\scriptsize$k-1$}] {};
\draw (k-1a) -- (L-1) -- (L);
\node (2k-1a) at (2.6,-1) {};
\node (2L-1) at (3,-1) [fill=black,draw,circle,label=below:{\scriptsize$2k-2$}] {};
\node (2L) at (4,-1) [fill=black,draw,circle,label=below:{\scriptsize$2k-1$}] {};
\draw [dotted] (2k-1a) -- (2L-1) -- (2L);
\draw (2k-1b) -- (2L-1) (L-1) -- (2L);
\draw [dotted] (k-1b) -- (L-1) (2L-1) -- (L); 
\node (2L+1) at (-1,-0) [fill=black,draw,circle,label=left:{\scriptsize$2k-1$}] {};
\node (2L+2) at (5,0) [fill=black,draw,circle,label=right:{\scriptsize$2k$}] {};
\draw [double,thick] (L) -- node[black,above right] {\scriptsize $\omega$} (2L+2) -- node[black,below right] {\scriptsize $-\omega$} (2L);
\draw [double,thick] (1) -- node[black,above left] {\scriptsize $\omega$} (2L+1) -- node[black,below left] {\scriptsize $\omega$} (L+1);
\end{tikzpicture}
\]
\[
\begin{tikzpicture}[scale=1.5]
\node (1) at (0,0) [fill=black,draw,circle,label=above:{\scriptsize$1$}] {};
\node (2) at (1,0) [fill=black,draw,circle,label=above:{\scriptsize$2$}] {};
\node (3) at (2,0) [fill=black,draw,circle,label=above:{\scriptsize$3$}] {};
\node (3a) at (2.4,0) {};
\node (3b) at (2.4,-0.4) {};
\draw (1) -- (2) -- (3) -- (3a); 
\draw (3) -- (3b);
\node (L+1) at (0,-1) [fill=black,draw,circle,label=below:{\scriptsize$k$}] {};
\node (L+2) at (1,-1) [fill=black,draw,circle,label=below:{\scriptsize$k+1$}] {};
\node (L+3) at (2,-1) [fill=black,draw,circle,label=below:{\scriptsize$k+2$}] {};
\node (k+3a) at (2.4,-1) {};
\node (k+3b) at (2.4,-0.6) {};
\draw [dotted] (L+1) -- (L+2) -- (L+3) -- (k+3a);
\draw [dotted] (L+3) -- (k+3b);
\draw (1) -- (L+2) (2) -- (L+3);
\draw [dotted] (L+1) -- (2) (L+2) -- (3);
\node (spacer) at (2.5,-0.5) {$\cdots$};
\node (k-1b) at (2.6,-0.4) {};
\node (2k-1b) at (2.6,-0.6) {};
\node (k-1a) at (2.6,0) {};
\node (L-1) at (3,0) [fill=black,draw,circle,label=above:{\scriptsize$k-2$}] {};
\node (L) at (4,0) [fill=black,draw,circle,label=above:{\scriptsize$k-1$}] {};
\draw (k-1a) -- (L-1) -- (L);
\node (2k-1a) at (2.6,-1) {};
\node (2L-1) at (3,-1) [fill=black,draw,circle,label=below:{\scriptsize$2k-2$}] {};
\node (2L) at (4,-1) [fill=black,draw,circle,label=below:{\scriptsize$2k-1$}] {};
\draw [dotted] (2k-1a) -- (2L-1) -- (2L);
\draw (2k-1b) -- (2L-1) (L-1) -- (2L);
\draw [dotted] (k-1b) -- (L-1) (2L-1) -- (L); 
\node (2L+1) at (-1,0) [fill=black,draw,circle,label=left:{\scriptsize$2k-1$}] {};
\node (2L+2) at (5,0) [fill=black,draw,circle,label=right:{\scriptsize$2k$}] {};
\draw [double,thick] (L) -- node[black,above right] {\scriptsize $\overline{\omega}$} (2L+2) -- node[black,below right] {\scriptsize $-\overline{\omega}$} (2L);
\draw [double,thick] (1) -- node[black,above left] {\scriptsize $\omega$} (2L+1) -- node[black,below left] {\scriptsize $\omega$} (L+1);
\end{tikzpicture}
\]
\caption{The families  of $2k$-vertex maximal cyclotomic $\mathcal{L}'$-graphs $T_{2k}^4$ and ($d=-7$ only) ${T_{2k}^4}'$, for $k\ge3$. \emph{(Where $\omega=\sqrt{-2},1/2+\sqrt{-7}/2$ for $d=-2,-7$ respectively.)}}
\label{T2k4figure}
\end{figure}

\begin{figure}[H]
\[
\begin{tikzpicture}[scale=1.5]
\node (1) at (0,0) [fill=white,inner sep=0pt,draw,circle,label=above:{\scriptsize$1$}] {$+$};
\node (2) at (1,0) [fill=black,draw,circle,label=above:{\scriptsize$2$}] {};
\node (3) at (2,0) [fill=black,draw,circle,label=above:{\scriptsize$3$}] {};
\node (3a) at (2.4,0) {};
\node (3b) at (2.4,-0.4) {};
\draw (1) -- (2) -- (3) -- (3a); 
\draw (3) -- (3b);
\node (L+1) at (0,-1) [fill=white,inner sep=0pt,draw,circle,label=below:{\scriptsize$k+1$}] {$+$};
\node (L+2) at (1,-1) [fill=black,draw,circle,label=below:{\scriptsize$k+2$}] {};
\node (L+3) at (2,-1) [fill=black,draw,circle,label=below:{\scriptsize$k+3$}] {};
\draw (1) -- (L+1);
\node (k+3a) at (2.4,-1) {};
\node (k+3b) at (2.4,-0.6) {};
\draw [dotted] (L+1) -- (L+2) -- (L+3) -- (k+3a);
\draw [dotted] (L+3) -- (k+3b);
\draw (1) -- (L+2) (2) -- (L+3);
\draw [dotted] (L+1) -- (2) (L+2) -- (3);
\node (spacer) at (2.5,-0.5) {$\cdots$};
\node (k-1b) at (2.6,-0.4) {};
\node (2k-1b) at (2.6,-0.6) {};
\node (k-1a) at (2.6,0) {};
\node (L-1) at (3,0) [fill=black,draw,circle,label=above:{\scriptsize$k-1$}] {};
\node (L) at (4,0) [fill=black,draw,circle,label=above:{\scriptsize$k$}] {};
\draw (k-1a) -- (L-1) -- (L);
\node (2k-1a) at (2.6,-1) {};
\node (2L-1) at (3,-1) [fill=black,draw,circle,label=below:{\scriptsize$2k-1$}] {};
\node (2L) at (4,-1) [fill=black,draw,circle,label=below:{\scriptsize$2k$}] {};
\draw [dotted] (2k-1a) -- (2L-1) -- (2L);
\draw (2k-1b) -- (2L-1) (L-1) -- (2L);
\draw [dotted] (k-1b) -- (L-1) (2L-1) -- (L); 
\node (2L+1) at (4.5,-0.5) [fill=black,draw,circle,label=right:{\scriptsize$2k+1$}] {};
\draw [double,thick] (L) -- node[black,above right] {\scriptsize $\omega$} (2L+1) -- node[black,below right] {\scriptsize $-\omega$} (2L);
\end{tikzpicture}
\]
\caption{The family of $2k+1$-vertex maximal cyclotomic $\mathcal{L}'$-graphs  $C_{2k}^{2+}$, for $k\ge1$. \emph{(Where $\omega=\sqrt{-2},1/2+\sqrt{-7}/2$ for $d=-2,-7$ respectively.)}}
\label{C2k+figure}
\end{figure}
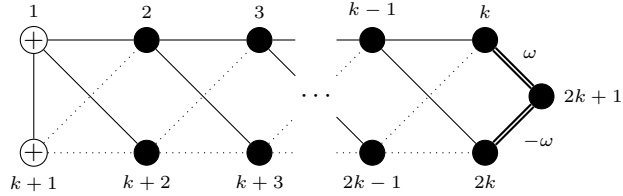

Summarising Theorem 1 of \cite{McSm1} and Theorems 8-11 of \cite{Ta2} we have:

\begin{theorem}\label{theclassification}
Let $R=\Z$ or $R=\mathcal{O}_{\Q(\sqrt{d})}$ for $d<0$, squarefree, $d\neq -1,-3$. If $A$ is a maximal cyclotomic Hermitian $R$-matrix, then $A$ has an $\mathcal{L}$-graph representative equivalent to one of the following: 
\begin{itemize}
\item The sporadic charged signed graph $S_7$, $S_8$, $S_8'$, $S_{14}$ or $S_{16}$ shown in Fig. \ref{sporadicCSGfigure};
\item A signed graph $T_{2k}$ for some $k\ge3$ as shown in Fig. \ref{T2kfigure};
\item A charged signed graph $C_{2k}^{++}$ or $C_{2k}^{+-}$ for some $k\ge2$ as shown in Fig. \ref{C2kfigure};
\item The sporadic $\mathcal{L}$-graph $S_2$ or ($d=-7,-15$ only) $S_2^*$ shown in Fig. \ref{weight4figure};
\item The sporadic $\mathcal{L}$-graph $S_2'$ ($d=-2,-11$ only)  or $S_4'$ ($d=-2,-11$ only)  shown in Fig. \ref{weight3figure}
\item The sporadic $\mathcal{L}'$-graph $S_4$ ($d=-2,-7$ only), $S_4^*$ ($d=-2$ only), $S_6^\dag$ ($d=-7$ only) or $S_8^*$ ($d=-2,-7$ only) shown in Fig. \ref{Sporadicw2figure};
\item An $\mathcal{L}'$-graph $T_{2k}^4$ ($d=-2,-7$ only) or ${T_{2k}^4}'$ ($d=-7$ only) for some $k\ge2$ as shown in Fig. \ref{T2k4figure};
\item An $\mathcal{L}'$-graph $C_{2k}^{2+}$ ($d=-2,-7$ only) for some $k\ge1$ as shown in Fig. \ref{C2k+figure}.
\end{itemize}
or is the matrix $(2)$ or $(-2)$.
\end{theorem}

\begin{figure}[H]
\[
\begin{tikzpicture}[scale=1.5]
\node (1) at (0,0) [fill=black,draw,circle,label=left:{$A$}] {};
\node (2) at (1,0) [fill=black,draw,circle] {};
\node (3) at (2,0) [fill=black,draw,circle] {};
\node (3a) at (2.4,0) {};
\node (3b) at (2.4,-0.4) {};
\draw (1) -- (2) -- (3) -- (3a); 
\draw (3) -- (3b);
\node (k+1) at (0,-1) [fill=black,draw,circle,label=left:{$B$}] {};
\node (k+2) at (1,-1) [fill=black,draw,circle] {};
\node (k+3) at (2,-1) [fill=black,draw,circle] {};
\node (k+3a) at (2.4,-1) {};
\node (k+3b) at (2.4,-0.6) {};
\draw [dotted] (k+1) -- (k+2) -- (k+3) -- (k+3a);
\draw [dotted] (k+3) -- (k+3b);
\draw [dotted] (k+2) -- (3);
\draw (1) -- (k+2);
\draw [dotted] (2) -- (k+1);
\draw (2) -- (k+3);
\node (spacer) at (2.5,-0.5) {$\cdots$};
\node (k-1b) at (2.6,-0.4) {};
\node (2k-1b) at (2.6,-0.6) {};
\node (k-1a) at (2.6,0) {};
\node (k-1) at (3,0) [fill=black,draw,circle] {};
\node (k) at (4,0) [fill=black,draw,circle] {};
\node (second1) at (5,0) [fill=black,draw,circle,label=right:{$A$}] {};
\draw (k-1a) -- (k-1) -- (k) -- (second1);
\node (2k-1a) at (2.6,-1) {};
\node (2k-1) at (3,-1) [fill=black,draw,circle] {};
\node (2k) at (4,-1) [fill=black,draw,circle] {};
\node (secondk+1) at (5,-1) [fill=black,draw,circle,label=right:{$B$}] {};
\draw [dotted] (2k-1a) -- (2k-1) -- (2k) -- (secondk+1);
\draw (2k-1b) -- (2k-1) (k-1) -- (2k) (k) -- (secondk+1);
\draw [dotted] (k-1b) -- (k-1) (2k-1) -- (k) (2k) -- (second1);
\end{tikzpicture}
\]
\caption{The family of $2k$-vertex maximal cyclotomic signed graphs $T_{2k}$, for $k\ge3$. \emph{(Where the two copies of vertices $A$ and $B$ should be identified.)}}
\label{T2kfigure}
\end{figure}
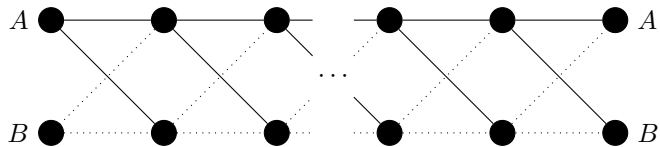

\begin{figure}[H]
\[
\begin{tikzpicture}[scale=1.5]
\node (1) at (0,0) [fill=white,inner sep=0pt,draw,circle] {$+$};
\node (2) at (1,0) [fill=black,draw,circle] {};
\node (3) at (2,0) [fill=black,draw,circle] {};
\node (3a) at (2.4,0) {};
\node (3b) at (2.4,-0.4) {};
\draw (1) -- (2) -- (3) -- (3a); 
\draw (3) -- (3b);
\node (L+1) at (0,-1) [fill=white,inner sep=0pt,draw,circle] {$+$};
\node (L+2) at (1,-1) [fill=black,draw,circle] {};
\node (L+3) at (2,-1) [fill=black,draw,circle] {};
\draw (1) -- (L+1);
\node (k+3a) at (2.4,-1) {};
\node (k+3b) at (2.4,-0.6) {};
\draw [dotted] (L+1) -- (L+2) -- (L+3) -- (k+3a);
\draw [dotted] (L+3) -- (k+3b);
\draw (1) -- (L+2) (2) -- (L+3);
\draw [dotted] (L+1) -- (2) (L+2) -- (3);
\node (spacer) at (2.5,-0.5) {$\cdots$};
\node (k-1b) at (2.6,-0.4) {};
\node (2k-1b) at (2.6,-0.6) {};
\node (k-1a) at (2.6,0) {};
\node (L-1) at (3,0) [fill=black,draw,circle] {};
\node (L) at (4,0) [fill=white,inner sep=0pt,draw,circle] {$+$};
\draw (k-1a) -- (L-1) -- (L);
\node (2k-1a) at (2.6,-1) {};
\node (2L-1) at (3,-1) [fill=black,draw,circle] {};
\node (2L) at (4,-1) [fill=white,inner sep=0pt,draw,circle] {$+$};
\draw [dotted] (2k-1a) -- (2L-1) -- (2L);
\draw (2k-1b) -- (2L-1) (L-1) -- (2L);
\draw [dotted] (k-1b) -- (L-1) (2L-1) -- (L); 
\draw [dotted] (L) -- (2L);
\node (C2k) at (-1,-0.5) {$C_{2k}^{++}:=$};
\end{tikzpicture}
\]
\[
\begin{tikzpicture}[scale=1.5]
\node (1) at (0,0) [fill=white,inner sep=0pt,draw,circle] {$+$};
\node (2) at (1,0) [fill=black,draw,circle] {};
\node (3) at (2,0) [fill=black,draw,circle] {};
\node (3a) at (2.4,0) {};
\node (3b) at (2.4,-0.4) {};
\draw (1) -- (2) -- (3) -- (3a); 
\draw (3) -- (3b);
\node (L+1) at (0,-1) [fill=white,inner sep=0pt,draw,circle] {$+$};
\node (L+2) at (1,-1) [fill=black,draw,circle] {};
\node (L+3) at (2,-1) [fill=black,draw,circle] {};
\draw (1) -- (L+1);
\node (k+3a) at (2.4,-1) {};
\node (k+3b) at (2.4,-0.6) {};
\draw [dotted] (L+1) -- (L+2) -- (L+3) -- (k+3a);
\draw [dotted] (L+3) -- (k+3b);
\draw (1) -- (L+2) (2) -- (L+3);
\draw [dotted] (L+1) -- (2) (L+2) -- (3);
\node (spacer) at (2.5,-0.5) {$\cdots$};
\node (k-1b) at (2.6,-0.4) {};
\node (2k-1b) at (2.6,-0.6) {};
\node (k-1a) at (2.6,0) {};
\node (L-1) at (3,0) [fill=black,draw,circle] {};
\node (L) at (4,0) [fill=white,inner sep=0pt,draw,circle] {$-$};
\draw (k-1a) -- (L-1) -- (L);
\node (2k-1a) at (2.6,-1) {};
\node (2L-1) at (3,-1) [fill=black,draw,circle] {};
\node (2L) at (4,-1) [fill=white,inner sep=0pt,draw,circle] {$-$};
\draw [dotted] (2k-1a) -- (2L-1) -- (2L);
\draw (2k-1b) -- (2L-1) (L-1) -- (2L);
\draw [dotted] (k-1b) -- (L-1) (2L-1) -- (L); 
\draw (L) -- (2L);
\node (C2k) at (-1,-0.5) {$C_{2k}^{+-}:=$};
\end{tikzpicture}
\]	
\caption{The families of $2k$-vertex maximal cyclotomic charged signed graphs $C_{2k}^{++}$ and $C_{2k}^{+-}$, for $k\ge2$.}
\label{C2kfigure}
\end{figure}
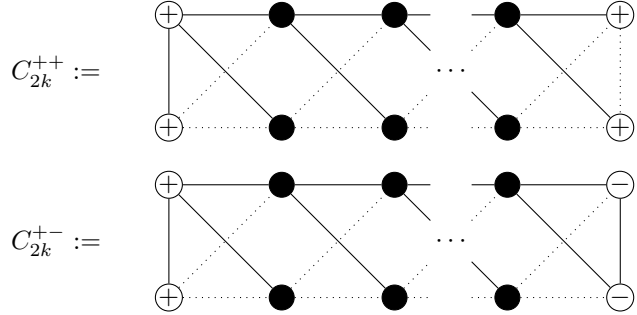

\begin{corollary}\label{wdmax4} Let $G$ be a cyclotomic $\mathcal{L}$-graph. Then each vertex of $G$ has weighted degree at most $4$.\end{corollary}

\section{Reduction to $\mathcal{L}'$-graphs}
Following \cite{Boyd} and \cite{Moss} we will consider a monic $f\in\Z[z]$ to have `small' Mahler measure if $M(f)<1.3$. A complete classification of integer symmetric matrices with small Mahler measure is given in \cite{McSm2}; we will ultimately prove Theorem \ref{lehmer271115} by showing that those are the only $R$-matrices with small Mahler measure. In this section we reduce the problem to finding $\mathcal{L}'$-graphs with small Mahler measure.

It will often be convenient to specify the edges of an $\mathcal{L}$-graph only up to their weight; we describe such a representation as the \emph{form} of the graph. Edges without an explicit label will be indicated by dashes ($\begin{tikzpicture}
\node (a) at (0,0) [fill=black,draw,shape=circle] {};
\node (b) at (0.9,0) [fill=black,draw,shape=circle] {};
\draw [dashed] (a) -- (b);
\end{tikzpicture}$, $\begin{tikzpicture}
\node (a) at (0,0) [fill=black,draw,shape=circle] {};
\node (b) at (0.9,0) [fill=black,draw,shape=circle] {};
\draw [dashed,double,thick] (a) -- (b);
\end{tikzpicture}$, $\begin{tikzpicture}
\node (a) at (0,0) [fill=black,draw,shape=circle] {};
\node (b) at (0.9,0) [fill=black,draw,shape=circle] {};
\draw [dashed] (a.0) -- (b.180) (a.15) -- (b.165) (a.345) -- (b.195);
\end{tikzpicture}$ for edges from $\mathcal{L}_1, \mathcal{L}_2, \mathcal{L}_3$ respectively) whilst an unspecified - possibly absent - edge will be shown as $\begin{tikzpicture}
\node (a) at (0,0) [fill=black,draw,shape=circle] {};
\node (b) at (1,0) [fill=black,draw,shape=circle] {};
\draw [dashdot] (a) -- (b); \end{tikzpicture}$. If a vertex is of unknown charge $c\in\{0,1,-1\}$ then we denote it by $\circledast$; a vertex known to be charged but of unknown polarity is denoted {\scriptsize\textcircled{$\pm$}}.

\begin{proposition}\label{diag2max} If $A$ is an $R$-matrix with a diagonal entry of modulus at least $3$, then $A$ has Mahler measure greater than $1.3$.\end{proposition}
\begin{proof} By interlacing $A$ contains a $1\times 1$ matrix $(n)$ for some $n\ge3$ which has Mahler measure $(n+\sqrt{n^2-4})/2\ge 2.618$.\end{proof} 

 \begin{proposition}\label{offdiag4max} If $A$ is an $R$-matrix with an off-diagonal entry of norm at least $5$, then $A$ has Mahler measure greater than $1.3$.\end{proposition}
 \begin{proof} By interlacing $A$ contains a $2\times 2$ matrix $A'=\left(\begin{array}{cc}b&a\\\overline{a}&c\end{array}\right)$ where $\mbox{norm}(a)=n\ge5$. By Proposition \ref{diag2max} we may assume $|b|,|c|$ at most $2$, but for each choice of $b,c$ the Mahler measure of $A'$ is increasing in $n$, with $n=5$ giving Mahler measure at least $2.36.$\end{proof}
 
 \begin{proposition}\label{diag2max2} If $A$ is a noncyclotomic $R$-matrix with a diagonal entry of modulus $2$, then $A$ has  Mahler measure greater than $1.3$.\end{proposition}
 \begin{proof} By interlacing, up to equivalence $A$, contains a $2\times 2$ matrix $A'=\left(\begin{array}{cc}2&a\\\overline{a}&b\end{array}\right)$ where by Propositions \ref{diag2max}, \ref{offdiag4max} we may assume $\mbox{norm}(a)=n\le4$ and $|b|\le2$. But any such choice of $n,b$ gives a matrix $A'$ with Mahler measure at least $1.722$.\end{proof}
 
 We may thus restrict our attention to $R$-matrices with an $\mathcal{L}$-graph representation.

 \begin{proposition}\label{w4edges}  Let $G$ be a noncyclotomic $\mathcal{L}$-graph with an edge of weight $4$. Then $G$ has Mahler measure greater than $1.3$.
 \end{proposition}
 \begin{proof} $G$ contains a subgraph $G'= \begin{tikzpicture}[scale=1]
\node (a) at (0,0) [fill=white,inner sep=0pt,draw,shape=circle] {$*$};
\node (b) at (1,0) [fill=white,inner sep=0pt,draw,shape=circle] {$*$};
\draw [double] (a.15) -- node[above] {\scriptsize $\omega$} (b.165);
\draw [double] (a.345) -- (b.195);\end{tikzpicture}$ for some $\omega\in\mathcal{L}_4$. If either vertex of $G'$ is charged then by Cor. \ref{wdmax4} it is noncyclotomic; any such graph has Mahler measure at least $2.08$. 
Otherwise $G'$ is uncharged and cyclotomic; by assumption it cannot be all of $G$, which must therefore induce a subgraph of the form 
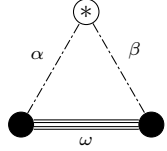
\[ \begin{tikzpicture}[scale=1]
\node (a) at (210:1) [fill=black,draw,shape=circle] {};
\node (b) at (330:1) [fill=black,draw,shape=circle] {};
\draw [double] (a.15) --  (b.165);
\draw [double] (a.345) -- node[below] {\scriptsize $\omega$} (b.195);
\node (c) at (90:1) [fill=white,draw,inner sep=1pt,shape=circle] {$*$};
\draw [dashdot] (a) -- node[above left] {\scriptsize $\alpha$} (c) -- node[above right] {\scriptsize $\beta$} (b);
\end{tikzpicture}\] where $\alpha,\beta\in\mathcal{L}$ are not both zero; any such graph also has Mahler measure at least $2.08$.\end{proof}
  
 We may thus restrict our attention to $\mathcal{L}$-graphs with all edges of weight at most $3$. 
  
 \begin{proposition}\label{w3edges} Let $G$ be a noncyclotomic $\mathcal{L}$-graph with an edge of weight $3$. Then $G$ has Mahler measure greater than $1.3$.\end{proposition}
 \begin{proof}
 By Theorem \ref{theclassification} any cyclotomic $\mathcal{L}$-graph containing a weight $3$ edge is (up to equivalence) $S_2'$, $S_4'$ or in an induced subgraph of $S_4'$. 
 If $G$ is a minimal noncyclotomic $\mathcal{L}$-graph containing $S_2'$ then it must be of the form
 
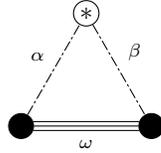
\[ \begin{tikzpicture}[scale=1]
\node (a) at (210:1) [fill=black,draw,shape=circle] {};
\node (b) at (330:1) [fill=black,draw,shape=circle] {};
\draw (a.0) -- (b.180) (a.340) -- node[below] {\scriptsize$\omega$} (b.200) (a.20) --  (b.160);
\node (c) at (90:1) [fill=white,draw,inner sep=1pt,shape=circle] {$*$};
\draw [dashdot] (a) -- node[above left] {\scriptsize $\alpha$} (c) -- node[above right] {\scriptsize $\beta$} (b);
\end{tikzpicture}\] where  $\alpha,\beta\in\mathcal{L}$ are not both zero; any such graph has Mahler measure at least $2.52$.
Otherwise, $G$ can have at most $5$ vertices; and (up to equivalence) is a subgraph of a $5$-vertex supergraph of $S_4'$. Given the constraints on possible edge labels this is a finite set, and we recover classes of $4,3$ and $2$-vertex minimal noncyclotomic $\mathcal{L}$-graphs, all having Mahler measure at least $1.56$.
 \end{proof}
 
We may thus restrict our attention to $\mathcal{L}'$-graphs, and hence $d=-2,-7$, since for other $d$ we have $\mathcal{L}'=\{0,\pm1,\pm2\}$ and so any $\mathcal{L}'$-graph is a charged signed graph.

\section{Small Minimal Noncyclotomic $\mathcal{L}'$-graphs with small Mahler measure}
In this section we will prove that for $d=-2,-7$:
\begin{proposition}\label{nosmallsmall} There are no minimal noncyclotomic $\mathcal{L}'$-graphs with at least one weight 2 edge label, at most ten vertices, and small ($<1.3$) Mahler measure. \end{proposition} 

\subsection{Growing Algorithms}

Given a maximal cyclotomic $\mathcal{L}$'-graph $G$ and an induced subgraph $G'$, we may recover $G$ from $G'$ by reintroducing the `missing' vertices one at a time, giving a sequence of cyclotomic supergraphs of $G'$ contained in $G$. We may thus recover all cyclotomic supergraphs of $G'$ by considering all possible sequences of additions of a new vertex to $G$. If an addition yields a connected cyclotomic graph we describe it as a \emph{cyclotomic addition}, otherwise as a \emph{noncyclotomic addition}: maximal graphs are therefore those which admit no cyclotomic addition.

Given an $n\times n$ $\mathcal{L}'$-matrix representative $A$ of $G$, the addition of an extra vertex is specified by a nonzero column vector $c\in {\mathcal{L}'}^n$ and a charge $x$ from charge set $X=\{0,1,-1\}$, giving a supermatrix
\[\left(\begin{array}{cccc}
A_{11} & \cdots & A_{1n} & c_1\\
\vdots & \ddots &\vdots & \vdots\\
A_{n1} & \cdots & A_{nn} & c_n\\
\overline{c_1} & \cdots & \overline{c_n} & x
\end{array}\right)\]

For a label set $L$ We define $C_n(L)$ to be the collection of nonzero vectors from $L^n$. The supermatrix formed from $A,c,x$ is equivalent to the one formed from $A,-c,x$ by switching at the extra vertex, so we restrict our attention to reduced column sets $C_{n'}(L)$ containing only one of each pair $\{c,-c\}$. Further, we may sometimes consider bounded column sets $C_{n'}^b(L)=\{ c\in C_{n'}(L)\,|\, \sum_1^n \mbox{Norm}(c_i) \le b\}$ when there are restrictions on the weighted degree of the vertex being added.

\subsection{Excluded Subgraphs}
\begin{figure}[H]

\[\mbox{$\begin{tikzpicture}
\node (b) at (330:1) [fill=black,draw,shape=circle] {};
\node (a) at (210:1) [fill=black,draw,shape=circle] {};
\node (c) at (90:1) [fill=black,draw,shape=circle] {};
\draw [double,dashed] (a) -- (c) -- (b);
\draw [dashed] (a) -- (b);
\node (label) at (0,-1) {$X_{3A}$};
\end{tikzpicture}$}
\hspace{3em}
\mbox{$\begin{tikzpicture}
\node (a) at (0,0.5) [fill=black,draw,shape=circle] {};
\node (b) at (1,0.5) [fill=black,draw,shape=circle] {};
\node (c) at (2,0.5) [fill=black,draw,shape=circle] {};
\node (d) at (3,0.5) [fill=white,inner sep=1.5pt,draw,shape=circle] {$*$};
\draw [double,dashed] (a) -- (b) -- (c);
\draw [dashed] (c) -- (d);
\node (label) at (1.5,-0.75) {$X_{4A}$};
\end{tikzpicture}$}\]
\caption{Excluded Subgraphs of type I}
\label{typeIgraphs}
\end{figure}
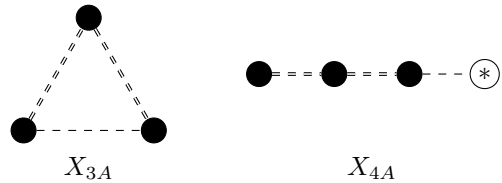

If an $\mathcal{L}'$-graph $G$ is minimal noncyclotomic then it cannot be a proper induced subgraph of any cyclotomic or minimal noncyclotomic $\mathcal{L}$'-graph; we describe $G$ as an excluded subgraph of type I. For $d=-2$, any $\mathcal{L}$'-graphs of the form $X_{3A}$, $X_{4A}$ given in Fig. \ref{typeIgraphs} are type I with Mahler measure at least $2.081\ldots$, whilst for $d=-7$ they are type I with Mahler measure at least $1.987\ldots$ and $2.081\ldots$ respectively.

\subsection{Excludable Subgraphs}
There exist cyclotomic $\mathcal{L}'$-graphs which are subgraphs of only finitely many maximal cyclotomic $\mathcal{L}'$-graphs; we describe these as type II graphs. If $H$ is type II, and the largest maximal cyclotomic $\mathcal{L}'$-graph to contain $H$ has $n$ vertices, then a minimal noncyclotomic $\mathcal{L}'$-graph containing $H$ has at most $n+1$ vertices. 

Any cyclotomic $\mathcal{L}'$-graph of one of the forms given in Fig. \ref{typeIIgraphs} is of type II:

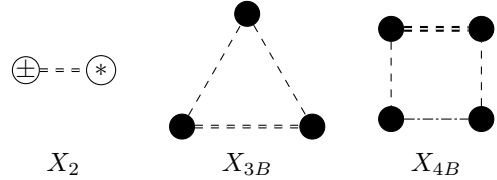
\begin{figure}[h]
\[
\mbox{$\begin{tikzpicture}
\node (a) at (0,0.5) [fill=white,inner sep=0pt,draw,shape=circle] {$\pm$};
\node (b) at (1,0.5) [fill=white,inner sep=1.5pt,draw,shape=circle] {$*$};
\draw [double,dashed] (a) -- (b);
\node (label) at (0.5,-0.75) {$X_2$};
\end{tikzpicture}$}\hspace{2em}
\mbox{$\begin{tikzpicture}
\node (b) at (330:1) [fill=black,draw,shape=circle] {};
\node (a) at (210:1) [fill=black,draw,shape=circle] {};
\node (c) at (90:1) [fill=black,draw,shape=circle] {};
\draw [dashed] (a) -- (c) -- (b);
\draw [dashed,double] (a) -- (b);
\node (label) at (0,-1) {$X_{3B}$};
\end{tikzpicture}$}\hspace{2em}
\mbox{$\begin{tikzpicture}[scale=1.2]
\node (b) at (1,0) [fill=black,draw,shape=circle] {};
\node (a) at (0,0) [fill=black,draw,shape=circle] {};
\node (c) at (0,-1) [fill=black,draw,shape=circle] {};
\node (d) at (1,-1) [fill=black,draw,shape=circle] {};
\draw [dashed] (a) -- (c) (d) -- (b);
\draw [dashed,double,thick] (a) -- (b);
\draw [dashdot] (c) -- (d);
\node (label) at (0.5,-1.5) {$X_{4B}$};
\end{tikzpicture}$}\] 
\caption{Excludable subgraphs of type II}
\label{typeIIgraphs}
\end{figure}

\begin{center}
\begin{tabular}{|c|c|c|}
\hline
Type II Graph & Maximal Supergraphs & Vertex Bound\\\hline
$X_2$ & $C_2^{2+}$, $S_4$ & 4\\
$X_{3B}$ & $S_4^*$ & 4\\
$X_4$ & $S_6^\dag$, $S_8^*$ & 8\\\hline
\end{tabular}
\end{center}

\subsection{The Search}
If $G$ is an $\mathcal{L}'$-graph with a weight 2 edge, then $G$ is equivalent to an $\mathcal{L}'$-graph inducing one of the $\mathcal{L}'$-graphs $H_i$ in Fig. \ref{theHi} as a subgraph.

\begin{figure}[H]
\[H_1:=\begin{tikzpicture}
\node (a) at (0,0) [fill=white,inner sep=0pt,draw,shape=circle] {$+$};
\node (b) at (1,0) [fill=white,inner sep=0pt,draw,shape=circle] {$+$};
\draw [double,thick] (a) --node[above] {\scriptsize $\omega$} (b);
\end{tikzpicture}\hspace{2em}
H_2:=\begin{tikzpicture}
\node (a) at (0,0) [fill=white,inner sep=0pt,draw,shape=circle] {$+$};
\node (b) at (1,0) [fill=white,inner sep=0pt,draw,shape=circle] {$-$};
\draw [double,thick] (a) --node[above] {\scriptsize $\omega$} (b);
\end{tikzpicture}\hspace{2em}
H_3:=\begin{tikzpicture}
\node (a) at (0,0) [fill=black,draw,shape=circle] {};
\node (b) at (1,0) [fill=black,draw,shape=circle] {};
\draw [double,thick] (a) --node[above] {\scriptsize $\omega$} (b);
\end{tikzpicture}\hspace{2em}
H_4:=\mbox{\begin{tikzpicture}
\node (a) at (0,0) [fill=white,inner sep=0pt,draw,shape=circle] {$+$};
\node (b) at (1,0) [fill=black,draw,shape=circle] {};
\draw [double,thick] (a) --node[above] {\scriptsize $\omega$} (b);
\end{tikzpicture}}
\]
\caption{Seed graphs $H_i$, where $\omega=\sqrt{-2}$ or $1/2+\sqrt{-7}/2$ for $d=-2,-7$ respectively.}
\label{theHi}
\end{figure}

We note that $H_1$ is noncyclotomic (with Mahler measure $1.883\ldots$) and clearly minimal. Let $G$ be a minimal noncyclotomic $\mathcal{L}'$-graph with vertices $x_1,\ldots,x_n$ for $3\le n\le 10$; then w.l.o.g. we may assume the subgraph induced on vertices $x_1,x_2$ is from the seed set $\{H_2,H_3,H_4\}$. By minimality, the subgraph on vertices $x_1,\ldots,x_{n-1}$ is cyclotomic. Thus $G$ can be recovered by a sequence of $n-3$ cyclotomic additions, followed by a noncyclotomic addition.

Starting with $\Sigma_2=\{H_2,H_3,H_4\}$, we take $\Sigma_i$ to be the set of all cyclotomic supergraphs of $\mathcal{L}'$-graphs in $\Sigma_{i-1}$, and $T_i$ to be their minimal noncyclotomic supergraphs; \emph{any} $n$-vertex minimal noncyclotomic $\mathcal{L}'$-graph with a weight 2 edge is then equivalent to some $\tau_n\in T_n$. This gives a large but finite search space; the sets $\Sigma_i, T_i$ can be generated more efficiently by avoiding the addition of vertices that would induce subgraphs of type I or II, or otherwise fail to be either cyclotomic or minimally noncyclotomic, as detailed below.

\begin{itemize}
\item[$n=3$] We grow the seed set $\Sigma_2$ by column set $C_{2'}(\mathcal{L'})$ and charge set $\{0,\pm1\}$. By brute force comparison of signed permutations, the sets $\Sigma_3$ and $T_3$ can be reduced modulo equivalence.

\item[$n=4$] We grow the sets $\Sigma_3$ by column set $C_{3'}(\mathcal{L'})$ and charge set $\{0,\pm1\}$; $\Sigma_4$ and $T_4$ can then be reduced modulo equivalence. 

\item[$n=5$] We grow the sets $\Sigma_4$ by column set $C_{4'}(\mathcal{L'})$ and charge set $\{0,\pm1\}$; $\Sigma_5$ and $T_5$ can then be reduced modulo equivalence. 

\item[$n=6$] To exclude $X_2$ we may restrict to $C_{5'}(\mathcal{L'})$, $X=\{0\}$ and $C_{5'}(\mathcal{L}_1\cup\{0\})$, $X=\{\pm1\}$. Brute-force equivalence testing is no longer feasible, but the sets $T_6$ are small enough to reduce by visual inspection of $\mathcal{L'}$-graphs.
 
\end{itemize}

For subsequent rounds, we note the following:

\begin{proposition}\label{bddgrow} By Cor. \ref{wdmax4}, any minimal noncyclotomic $\mathcal{L}'$-graph with seven or more vertices has all vertices of weighted degree at most four.\end{proposition}

\begin{itemize}
\item[$n=7,8,9$] To exclude $X_2$ and by proposition \ref{bddgrow} we may restrict to $C_{(n-1)'}^4(\mathcal{L}')$, $X=\{0\}$ and $C_{(n-1)'}^3(\mathcal{L}_1\cup\{0\})$, $X=\{\pm1\}$, and consider only supergraphs with vertices of weight at most four. In each round, and for both $d$, we find only a single minimal noncyclotomic supergraph.

\item[$n=10$] To exclude $X_2$ and by proposition \ref{bddgrow} we grow by $C_{9'}^3(\mathcal{L}_1\cup\{0\})$, $X=\{\pm1\}$; $C_{9'}^4(\mathcal{L}_1\cup\{0\})$, $X=\{0\}$; and $C_{9'}^4(\mathcal{L}_2\cup\{0\})$, $X=\{0\}$, considering only supergraphs with vertices of weight at most four. 

This leaves only addition of vectors from $C_{9'}^4(\mathcal{L}')$ with entries from both $\mathcal{L}_1$ and $\mathcal{L}_2$, the majority of which can be discarded as they would induce a subgraph of form $X_2$, $X_{4A}$ or $X_{4B}$. 

For both $d$, we find only a single minimal noncyclotomic supergraph.
\end{itemize}

\subsection{Results}
The results of the above search are summarised in Fig. \ref{searchresults}; for each $i$, the least Mahler measure listed was attained by $\mathcal{L}'$-graphs for both $d=-2$ and $d=-7$. For a full list of representatives, see Chapter 6 of \cite{Ta1}. Since no $\mathcal{L}'$-graph found had Mahler measure less than $1.35$, and the excluded graphs of type I had Mahler measure at least $1.987\ldots$, Proposition \ref{nosmallsmall} holds.

\begin{figure}[H]
\begin{center}
\begin{tabular}{|c|c|c|c|}
\hline
$i$ & \multicolumn{2}{|c|}{$|T_i|$} & $\displaystyle\min_{A\in T_i}{M(A)}$ \\
 & $d=-2$ & $d=-7$ & \\\hline
3 & 34 & 67 & 1.401\ldots\\
4 & 51 & 61 & 1.401\ldots\\
5 & 14 & 25 & 1.351\ldots\\
6 & 12 & 17 & 1.401\ldots\\
7 & 1 & 1 & 1.506\ldots\\
8 & 1 & 1 & 1.458\ldots\\
9 & 1 & 1 & 1.425\ldots\\
10 & 1 & 1 & 1.401\ldots\\\hline
\end{tabular}
\end{center}
\caption{Least Mahler measures of small $\mathcal{L}'$-graphs with at least one weight 2 edge.}
\label{searchresults}
\end{figure}

\section{Large Minimal Noncyclotomic $\mathcal{L}'$-graphs with small Mahler measure}
In this section we will prove that for $d=-2,-7$,

\begin{proposition}\label{v10max} If $G$ is a minimal noncyclotomic $\mathcal{L}'$-graph, then $G$ has at most ten vertices.\end{proposition}

\subsection{Supersporadics}\label{supersporadicsection}

For $d=-2,-7$, let $\mathcal{S}_d$ be the set of sporadic maximal cyclotomic graphs with edges of weight at most $2$; from Theorem \ref{theclassification} we have 

\[\mathcal{S}_{-2}=\{S_4,S_4^*,S_7,S_8,S'_8,S_8^*,S_{14},S_{16}\},\,\, \mathcal{S}_{-7}=\{S_4,S^\dag_6,S_7,S_8,S'_8,S^*_8,S_{14}, S_{16}\}\]

We then describe an $n$-vertex minimal noncyclotomic $\mathcal{L}'$-graph as \emph{supersporadic} if it has a connected subgraph with $n-1$ vertices which is equivalent to a subgraph $H$ of some $G\in\mathcal{S}_d$.

The set of supersporadic minimal noncyclotomic $\mathcal{L}'$-graphs is finite, and could (in principle) be computed from the set of all subgraphs of each $G\in\mathcal{S}_d$ by considering all possible single-vertex additions to each such subgraph. By restricting to $G$ with at least $11$ vertices, we need only consider connected minimal noncyclotomic $\mathcal{L}'$-graphs $G$ obtained by the noncyclotomic addition of a vertex $x$ to a $k$-vertex subgraph $H$ of $S_{14}$ or $S_{16}$, for $k\ge10$. 

By Proposition \ref{bddgrow}, $x$ has weighted degree at most $4$, so we may restrict to addition vectors $c$ from $\mathcal{C}_{k'}^4(\mathcal{L}')$. Proposition \ref{v10max} holds for charged signed graphs by the results of \cite{McSm2}, so we need only identify $\mathcal{L}'$-graphs with at least one edge of weight $2$: such an edge must be incident at $x$ since all edges in $H$ have weight $1$. Thus to exclude subgraphs of the form $X_2$, $x$ must be uncharged. 

If $H$ is connected, then a subgraph of form $X_{3B}$ or $X_{4B}$ will be induced unless $c\in C_{k'}^4(\mathcal{L}_2\cup\{0\})$, whilst if $H$ is disconnected then any supergraph will induce a subgraph of form $X_{4B}$ unless one of the connected components of $H$ is a singleton vertex. 

For each $10\le k \le 16$, it is then feasible to grow supergraphs of the connected representatives of the $k$-vertex subgraphs of $S_{14}$ and $S_{16}$ by column set $C_{k'}^4(\mathcal{L}_2\cup\{0\})$ with $X=\{0\}$, and the suitable disconnected representatives (of which there are very few) by column set $C_{k'}^4(\mathcal{L}')$ with $X=\{0\}$. This process confirms that no such supergraph is minimal noncyclotomic, and so Proposition \ref{v10max} holds for supersporadic $G$.

\subsection{Non-supersporadics}
Let $G$ be an $n\ge 11$-vertex minimal noncyclotomic $\mathcal{L}$-graph with a weight 2 edge. By minimality, each of the $(n-1)$-vertex subgraphs $G'_i$ of $G$ must be cyclotomic. Since we have shown that $G$ is not supersporadic, by Theorem \ref{theclassification} each of the $G'_i$ (and hence their subgraphs) are equivalent to subgraphs of some $T_{2k}$, $C^{+\pm}_{2k}$, $C^{2+}_{2k}$, $T^4_{2k}$ or ${T^4_{2k}}'$. We will prove the following:

\begin{proposition}\label{newprop8} Let $G$ be an $(n\ge 11)$-vertex connected $\mathcal{L}'$-graph such that every proper connected subgraph of $G$ is equivalent to a subgraph of some $T_{2k}$, $C^{+\pm}_{2k}$, $C^{2+}_{2k}$, $T^4_{2k}$ or ${T^4_{2k}}'$. Then $G$ is also equivalent to a subgraph of some $T_{2k}$, $C^{+\pm}_{2k}$, $C^{2+}_{2k}$, $T^4_{2k}$ or ${T^4_{2k}}'$.
\end{proposition}

The proof of this result in \cite{McSm2} for connected charged signed graphs mostly generalises in a straightforward way to $\mathcal{L}'$-graphs; we note the necessary changes.

\subsubsection{Profiles and Ranks}

\begin{proposition}\label{T2k4profile} The $2k$-vertex graphs $T^4_{2k},{T^4_{2k}}'$ have profiles of rank $k+1$:
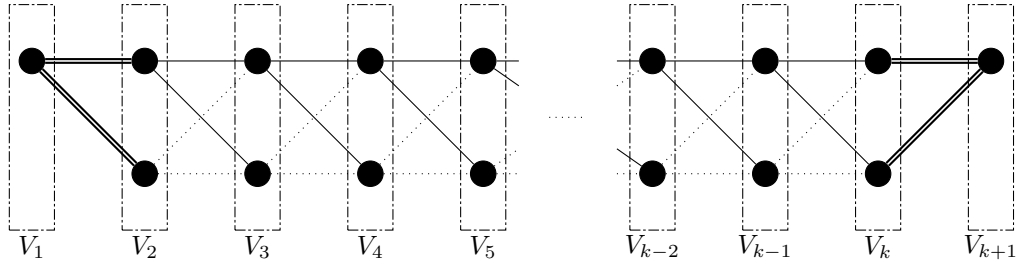
\[
\begin{tikzpicture}[scale=1.5]
\node (2l+1) at (-1,1) [fill=black,draw,shape=circle] {} ;
\node (k+1) at (0,0) [fill=black,draw,shape=circle] {} ;
\node (1) at (0,1) [fill=black,draw,shape=circle] {};
\draw [double,thick] (2l+1) -- (1);
\draw [double,thick] (2l+1) -- (k+1);
\node (k+2) at (1,0) [fill=black,draw,shape=circle] {} ;
\node (2) at (1,1) [fill=black,draw,shape=circle] {} ;
\node (k+3) at (2,0) [fill=black,draw,shape=circle] {} ;
\node (3) at (2,1) [fill=black,draw,shape=circle] {} ;
\node (k+4) at (3,0) [fill=black,draw,shape=circle] {} ;
\node (4) at (3,1) [fill=black,draw,shape=circle] {} ;
\node (k+4') at (3.4,0) {} ;
\node (4') at (3.4,1) {} ;
\node (x) at (3.4,0.72) {};
\node (y) at (3.4,0.28) {};
\draw [dotted] (k+1) -- (k+2) -- (k+3) -- (k+4) -- (k+4');
\draw (1) -- (2) -- (3) -- (4) -- (4');
\draw (1) -- (k+2) (2) -- (k+3) (3) -- (k+4) (4) -- (x);
\draw [dotted] (k+1) -- (2) (k+2) -- (3) (k+3) -- (4) (k+4) -- (y);
\node (a) at (3.5,0.5) {};
\node (b) at (4,0.5) {};
\draw [dotted] (a) -- (b);

\node (x') at (4.1,0.72) {};
\node (y') at (4.1,0.28) {};
\node (5') at (4.1,1) {};
\node (k+5') at (4.1,0) {};
\node (k+5) at (4.5,0) [fill=black,draw,shape=circle] {} ;
\node (5) at (4.5,1) [fill=black,draw,shape=circle] {} ;
\draw [dotted] (x') -- (5) (k+5') -- (k+5);
\draw (y') -- (k+5) (5') -- (5);
\node (6) at (5.5,1) [fill=black,draw,shape=circle] {};
\node (7) at (6.5,1) [fill=black,draw,shape=circle] {};

\node (k+6) at (5.5,0) [fill=black,draw,shape=circle] {};
\node (k+7) at (6.5,0) [fill=black,draw,shape=circle] {};

\draw [dotted] (k+5) -- (k+6) -- (k+7);
\draw (5) -- (6) -- (7);
\draw (5) -- (k+6) (6) -- (k+7) (7);
\draw [dotted] (k+5) -- (6) (k+6) -- (7) (k+7);
\node (2l+2) at (7.5,1) [fill=black,draw,shape=circle] {};
\draw [double,thick] (2l+2) -- (7);
\draw [double,thick] (2l+2) -- (k+7);

\foreach \x in {-1cm,0cm,1cm,2cm,3cm,4.5cm,5.5cm,6.5cm,7.5cm}
	\draw [dashdot] (-0.2cm+\x,-0.5cm) rectangle (0.2cm+\x,1.5cm);

\node (V0) at (-1,-0.65) {$V_1$};
\node (V1) at (0,-0.65) {$V_2$};
\node (V2) at (1,-0.65) {$V_3$};
\node (V3) at (2,-0.65) {$V_4$};
\node (V4) at (3,-0.65) {$V_5$};
\node (V5) at (4.5,-0.65) {$V_{k-2}$};
\node (V6) at (5.5,-0.65) {$V_{k-1}$};
\node (V7) at (6.5,-0.65) {$V_{k}$};
\node (V8) at (7.5,-0.65) {$V_{k+1}$};
\end{tikzpicture}
\]
\end{proposition}

\begin{proposition}\label{C2k2profile} The $2k+1$-vertex graphs $C_{2k}^{2+}$ have profiles of rank $k+1$:
\[
\begin{tikzpicture}[scale=1.5]
\node (k+1) at (0,0) [fill=white,inner sep=0pt,draw,shape=circle] {$+$} ;
\node (1) at (0,1) [fill=white,draw,inner sep=0pt,shape=circle] {$+$};
\node (k+2) at (1,0) [fill=black,draw,shape=circle] {} ;
\node (2) at (1,1) [fill=black,draw,shape=circle] {} ;
\node (k+3) at (2,0) [fill=black,draw,shape=circle] {} ;
\node (3) at (2,1) [fill=black,draw,shape=circle] {} ;
\node (k+4) at (3,0) [fill=black,draw,shape=circle] {} ;
\node (4) at (3,1) [fill=black,draw,shape=circle] {} ;
\node (k+4') at (3.4,0) {} ;
\node (4') at (3.4,1) {} ;
\node (x) at (3.4,0.72) {};
\node (y) at (3.4,0.28) {};
\draw [dotted] (k+1) -- (k+2) -- (k+3) -- (k+4) -- (k+4');
\draw (1) -- (2) -- (3) -- (4) -- (4') (1) -- (k+1);
\draw (1) -- (k+2) (2) -- (k+3) (3) -- (k+4) (4) -- (x);
\draw [dotted] (k+1) -- (2) (k+2) -- (3) (k+3) -- (4) (k+4) -- (y);
\node (a) at (3.5,0.5) {};
\node (b) at (4,0.5) {};
\draw [dotted] (a) -- (b);

\node (x') at (4.1,0.72) {};
\node (y') at (4.1,0.28) {};
\node (5') at (4.1,1) {};
\node (k+5') at (4.1,0) {};
\node (k+5) at (4.5,0) [fill=black,draw,shape=circle] {} ;
\node (5) at (4.5,1) [fill=black,draw,shape=circle] {} ;
\draw [dotted] (x') -- (5) (k+5') -- (k+5);
\draw (y') -- (k+5) (5') -- (5);
\node (6) at (5.5,1) [fill=black,draw,shape=circle] {};
\node (7) at (6.5,1) [fill=black,draw,shape=circle] {};
\node (k+6) at (5.5,0) [fill=black,draw,shape=circle] {};
\node (k+7) at (6.5,0) [fill=black,draw,shape=circle] {};
\node (2l+2) at (7.5,1) [fill=black,draw,shape=circle] {};
\draw [double,thick] (2l+2) -- (7);
\draw [double,thick] (2l+2) -- (k+7);
\draw [dotted] (k+5) -- (k+6) -- (k+7);
\draw (5) -- (6) -- (7);
\draw (5) -- (k+6) (6) -- (k+7) (7);
\draw [dotted] (k+5) -- (6) (k+6) -- (7) (k+7);

\foreach \x in {0cm,1cm,2cm,3cm,4.5cm,5.5cm,6.5cm,7.5cm}
	\draw [dashdot] (-0.2cm+\x,-0.5cm) rectangle (0.2cm+\x,1.5cm);

\node (V1) at (0,-0.65) {$V_1$};
\node (V2) at (1,-0.65) {$V_2$};
\node (V3) at (2,-0.65) {$V_3$};
\node (V4) at (3,-0.65) {$V_4$};
\node (V5) at (4.5,-0.65) {$V_{k-2}$};
\node (V6) at (5.5,-0.65) {$V_{k-1}$};
\node (V7) at (6.5,-0.65) {$V_{k}$};
\node (V8) at (7.5,-0.65) {$V_{k+1}$};

\end{tikzpicture}
\]
\end{proposition}

Lemma 6 in \cite{McSm2} then generalises to 

\begin{lemma}\label{uniqueprofile} Let $G$ be equivalent to a connected subgraph of one of $T_{2k},C^{+\pm}_{2k},C^{2+}_{2k},T^4_{2k}$ or ${T^4_{2k}}'$. If $G$ has path rank at least 5 then this equals its profile rank, and its columns are uniquely determined. Moreover, their order is determined up to reversal or cycling.
\end{lemma}

with the proof for $C^{2+}_{2k}$, $T^4_{2k}$ and ${T^4_{2k}}'$ proceeding as for $T_{2k}$ and $C^{+\pm}_{2k}$.

\subsubsection{Subgraph Conditions}

Any induced 4-cycle in a subgraph of rank at least 5 must then be one of the following:
\paragraph{Hourglass 4-cycles} Underlying graph of form 
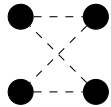
\[\begin{tikzpicture}
\node (a) at (0,0) [fill=black,draw,shape=circle] {};
\node (b) at (1,0) [fill=black,draw,shape=circle] {};
\node (c) at (0,-1) [fill=black,draw,shape=circle] {};
\node (d) at (1,-1) [fill=black,draw,shape=circle] {};
\draw [dashed] (a) -- (b) -- (c) -- (d) -- (a);
\end{tikzpicture}
\]
\paragraph{Parallelogram 4-cycles} Underlying graph of one of the forms 
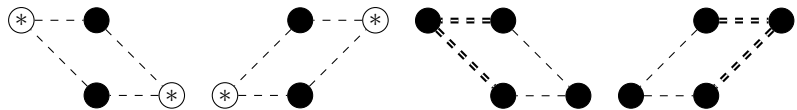
\[
\mbox{\begin{tikzpicture}
\node (a) at (0,0) [fill=white,inner sep=1pt,draw,shape=circle] {$*$};
\node (b) at (1,0) [fill=black,draw,shape=circle] {};
\node (c) at (1,-1) [fill=black,draw,shape=circle] {};
\node (d) at (2,-1) [fill=white,inner sep=1pt,draw,shape=circle] {$*$};
\draw [dashed] (a) -- (b) -- (d) -- (c) -- (a);
\end{tikzpicture}}
\hspace{1em}
\mbox{\begin{tikzpicture}
\node (a) at (0,0) [fill=white,inner sep=1pt,draw,shape=circle] {$*$};
\node (b) at (-1,0) [fill=black,draw,shape=circle] {};
\node (c) at (-1,-1) [fill=black,draw,shape=circle] {};
\node (d) at (-2,-1) [fill=white,inner sep=1pt,draw,shape=circle] {$*$};
\draw [dashed] (a) -- (b) -- (d) -- (c) -- (a);
\end{tikzpicture}}
\hspace{1em}
\mbox{\begin{tikzpicture}
\node (a) at (0,0) [fill=black,draw,shape=circle] {};
\node (b) at (1,0) [fill=black,draw,shape=circle] {};
\node (c) at (1,-1) [fill=black,draw,shape=circle] {};
\node (d) at (2,-1) [fill=black,draw,shape=circle] {};
\draw [dashed] (b) -- (d) -- (c);
\draw [dashed,double,thick] (b) -- (a) -- (c);
\end{tikzpicture}}
\hspace{1em}
\mbox{\begin{tikzpicture}
\node (a) at (0,0) [fill=black,draw,shape=circle] {};
\node (b) at (-1,0) [fill=black,draw,shape=circle] {};
\node (c) at (-1,-1) [fill=black,draw,shape=circle] {};
\node (d) at (-2,-1) [fill=black,draw,shape=circle] {};
\draw [dashed] (b) -- (d) -- (c);
\draw [dashed,double,thick] (b) -- (a) -- (c);
\end{tikzpicture}}
\]

\paragraph{Triangular 4-cycles} Underlying graph of one of the forms
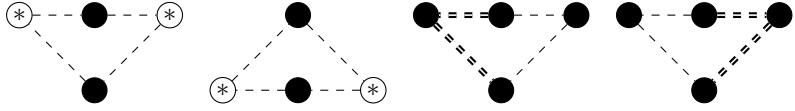
\[
\mbox{\begin{tikzpicture}
\node (a) at (0,0) [fill=white,inner sep=1pt,draw,shape=circle] {$*$};
\node (b) at (1,0) [fill=black,draw,shape=circle] {};
\node (c) at (1,-1) [fill=black,draw,shape=circle] {};
\node (d) at (2,0) [fill=white,inner sep=1pt,draw,shape=circle] {$*$};
\draw [dashed] (a) -- (b) -- (d) -- (c) -- (a);
\end{tikzpicture}}
\hspace{1em}
\mbox{\begin{tikzpicture}
\node (a) at (0,0) [fill=white,inner sep=1pt,draw,shape=circle] {$*$};
\node (b) at (1,1) [fill=black,draw,shape=circle] {};
\node (c) at (1,0) [fill=black,draw,shape=circle] {};
\node (d) at (2,0) [fill=white,inner sep=1pt,draw,shape=circle] {$*$};
\draw [dashed] (a) -- (b) -- (d) -- (c) -- (a);
\end{tikzpicture}}
\hspace{1em}
\mbox{\begin{tikzpicture}
\node (a) at (0,0) [fill=black,draw,shape=circle] {};
\node (b) at (1,0) [fill=black,draw,shape=circle] {};
\node (c) at (1,-1) [fill=black,draw,shape=circle] {};
\node (d) at (2,0) [fill=black,draw,shape=circle] {};
\draw [dashed] (b) -- (d) -- (c);
\draw [dashed,double,thick] (b) -- (a) -- (c);
\end{tikzpicture}}
\hspace{1em}
\mbox{\begin{tikzpicture}
\node (a) at (0,0) [fill=black,draw,shape=circle] {};
\node (b) at (-1,0) [fill=black,draw,shape=circle] {};
\node (c) at (-1,-1) [fill=black,draw,shape=circle] {};
\node (d) at (-2,0) [fill=black,draw,shape=circle] {};
\draw [dashed] (b) -- (d) -- (c);
\draw [dashed,double,thick] (b) -- (a) -- (c);
\end{tikzpicture}}
\]

We fix the numbering of vertices of $T_{2k}^4$, ${T_{2k}^4}'$, $C_{2k}^{2+}$ and their subgraphs as in Figures \ref{T2k4figure} and \ref{C2k+figure}. For $d=-2$, we define an edge to be \emph{positive} if it has label $+1$ or $\omega=\sqrt{-2}$; otherwise (label from $\{-1,-\sqrt{-2}\}$) we call it \emph{negative}. For $d=-7$, we define an edge to be \emph{positive} if it has label from $\{+1,\omega=1/2+\sqrt{-7}/2,\overline{\omega}\}$ or \emph{negative} if it has label from $\{-1,-\omega,-\overline{\omega}\}$.

Proposition 7 of \cite{McSm2} can then be extended to include the following cases:

\begin{proposition}\label{subgraphprop}
(iii) Let $H$ be an uncharged $\mathcal{L}$-graph of rank at least 5 that has, for some $k$, an underlying graph of the same form as a subgraph of $T_{2k}^4$ or $T_{2k}^{4'}$, as drawn in Proposition \ref{T2k4profile}. Then $H$ is equivalent to a subgraph $G$ of $T_{2k}^4$ or $T_{2k}^{4'}$ if and only if
\begin{itemize}
\item The hourglass 4-cycles all have an even number of positive edges;
\item The parallelogram 4-cycles all have an odd number of positive edges;
\item The triangular 4-cycles all have an odd number of positive edges.
\end{itemize}
(iv) Let $H$ be a charged $\mathcal{L}$-graph of rank at least 5 that has, for some $k$, an underlying graph of the same form as a subgraph of $C_{2k}^{2+}$ as drawn in Proposition \ref{C2k2profile}. Then $H$ is equivalent to a subgraph $G$ of $C_{2k}^{2+}$ if and only if 
\begin{itemize}
\item The hourglass 4-cycles all have an even number of positive edges;
\item The parallelogram 4-cycles all have an off number of positive edges;
\item The triangular 4-cycles all have an odd number of positive edges;
\item The triangles containing two charged vertices in the subgraph have the property that if the charges are positive (respectively negative) then the triangle has an even number of positive (resp. negative) edges.
\end{itemize}
\end{proposition}

\begin{proof}
We first show that the conditions given in Proposition \ref{subgraphprop} are necessary. Since $H$ has rank at least 5, by Lemma \ref{uniqueprofile} the columns of its profile are uniquely determined. Thus as drawn in Propositions \ref{T2k4profile}-\ref{C2k2profile} each 4-cycle of $H$ is either
\begin{itemize}
\item an hourglass
\end{itemize}
or
\begin{itemize}
\item a parallelogram 4-cycle or triangular 4-cycle. (Interchanging the position of conjugate vertices in the drawing may cause parallelograms to become triangular, and vice versa).
\end{itemize}
Since each 4-cycle is even length and contains zero or two edges of weight 2, the equivalence relation operations (permutation, switching, conjugation) will preserve the parity of the number of positive edges in each cycle, proving necessity.
We now assume that the given conditions hold, and prove that they are sufficient: that our given subgraph is equivalent to a subgraph of $T_{2k},C^{+\pm}_{2k},C^{2+}_{2k},T^4_{2k}$ or ${T^4_{2k}}'$. 
To do this, we need to embed an $\mathcal{L}'$-graph equivalent to $H$ into one of $T_{2k},C^{+\pm}_{2k},C^{2+}_{2k},T^4_{2k}$ or ${T^4_{2k}}'$ so that the resultant embedding $G$ inherits its edge and vertex signs from the $\mathcal{L}$-graph it is embedded into. 
Cases (i) and (ii) hold by Proposition 7 of \cite{McSm2}; for (iii) and (iv) we may assume that $H$ contains at least one edge of weight 2 else the conditions for that result are met with $H$ equivalent to a subgraph of $T_{2k},C_{2k}^{++}$ or $C_{2k}^{+-}$.
\begin{itemize}
\item[(iii)] Given that $H$ contains a weight 2 edge it cannot be equivalent to a a subgraph of $T_{2k},C_{2k}^{++}$ or $C_{2k}^{+-}$, and as it is uncharged we therefore seek to embed an equivalent of $H$ in $T_{2k}^4$ (or ${T_{2k}^4}'$ for $d=-7$).
 
Let $P$ be a maximum-length chordless path or cycle in $H$; since no chordless cycle in the underlying graph of $T_{2k}^4$ or $T_{2k}^{4'}$ has length greater than 4 but $H$ has rank at least 5, $P$ is necessarily a chordless path. Let it have length $l'$, joining vertices $v_1,\ldots,v_l$; by switching, we can ensure that it has all edge labels positive.

Let $e$ be an edge of weight 2 in $H$; w.l.o.g we may draw $H$ such that $e$ is the leftmost edge joining vertices 1 and $2k-1$ (as numbered in Figure \ref{T2k4figure}). Any longest rational integer path $P'$ in $T_{2k}^4$ or $T_{2k}^{4'}$ is at most $k-1$ vertices long. Consider its leftmost vertices $v_1,v_2$. If $v_1=1$ then a longer chordless path is obtained by starting at $2k-1$ then proceeding as in $P'$ via 1; if column $V_1\neq\{1\}$ then both vertices $1$ and $k$ are in $H$ (else redraw and take $k$ as $1$) so there is a longer path through $2k-1,1,v_2$ then proceeding as in $P'$. So the longest chordless path cannot have all edges rational integers and we may assume that the first edge of $P$ is of weight 2. 

Now either the edge between vertices $v_{l-1},v_l$ of $P$ is weight 2, or it isn't. If it is, we may embed $P$ into the top edge of $T=T_{2l-2}^4$ (or, for $d=-7$, $T={T_{2l-2}^4}'$ if the second weight 2 edge label is complex conjugate to the first); otherwise, embed into $T=T_{2l}^4$. In either case, all the relevant edges are positive as required. We may proceed as in case (i) in \cite{McSm2}; the next two paragraphs are essentially identical to that proof.

We can now embed into $T$ those conjugates of $v_1,\ldots,v_l$ that are present in $H$, by placing them in their appropriate columns on the bottom row of $T$: note that triangular 4-cycles in $H$ may become parallelogram 4-cycles, and vice versa, by this process (if $P$ moved between the top and bottom rows of the original drawing). This induces an embedding $G$ of $H$ in $T$, though without the signs of the edges yet agreeing. To achieve this agreement, we switch at these newly embedded vertices, if necessary, to ensure that all edges of negative slope have positive sign. We also switch at any vertex in the bottom row that has no incident edge of negative slope, if necessary, to ensure that the incident edge of positive slope has negative sign. 

We next claim that, after making these switchings, all edges of the embedding $G$ do indeed have the same sign as the edges of $T$. First consider an edge of $G$ of positive slope. If not already made to have negative sign, such an edge must be part of a triangular 4-cycle where the two horizontal edges and the edge of negative slope all have positive sign. Hence, by the stated triangular 4-cycle condition, the edge of positive slope must have negative sign. (Note that because both the stated parallelogram 4-cycle condition and the triangular 4-cycle condition hold for $H$, the triangular 4-cycle condition holds for $G$.) Finally, every horizontal edge on the second row is part of an hourglass 4-cycle, which implies that it must have negative sign.

\item[(iv)] Again, consider a maximum length chordless path $P$ in $H$. If no vertex is charged then $P$ could be embedded in $T_{2k}^4$ or $T_{2k}^{4'}$. So we may assume that $P$ contains a charged vertex: by the profile of $C_{2k}^{2\pm}$, this must be an end vertex of $P$. Further, by maximality, $P$ must terminate with a weight 2 edge. Negating if necessary we may assume that the charged vertex is positive, and by switching we may ensure that all edges of $P$ are positive, and by taking the complex conjugate if necessary that the weight 2 edge is $\omega$. Then such a $P$ with $k'$ vertices can be embedded sign-consistently into the top row of $C_{2(k'-1)}^{2+}$. We then proceed as in (iii), which ensures that all horizontal edges, and those of positive or negative slope, have the correct sign. Finally, the triangle condition ensures that the vertical edge must have positive sign as required. 
\end{itemize}
\end{proof}

To complete the proof of Proposition \ref{newprop8} we need only consider $\mathcal{L}'$-graphs with at least one weight $2$ edge. Let $G$ be such a graph. 

\begin{proposition}\label{no5cycle}  For $n\ge5$, $G$ cannot contain a chordless $n$-cycle. \end{proposition}

\begin{proof}
Let $G$ contain a chordless $n$-cycle on vertices $v_1,\ldots,v_n$. Further, by assumption there exist vertices $v,v'$ (possibly in $\{v_1,\ldots,v_n\}$) such that $e_{v,v'}\in\mathcal{L}_2$. Now let $G'$ be the smallest connected subgraph of $G$ to include all of $v_1,\ldots,v_n,v,v'$. If $G'$ is a proper subgraph of $G$, then we have a contradiction: $G'$ must be equivalent to a subgraph of some $T_{2k},C^{+\pm}_{2k},C^{2+}_{2k},T^4_{2k}$ or $T^{4'}_{2k}$, but none of those contain both an $\mathcal{L}_2$ edge and a chordless $n$-cycle on more than 4 vertices. Thus $G'=G$, and deleting any vertex not from $\{v_1,\ldots,v_n,v,v'\}$ gives a disconnected graph.

If $v,v'\in\{v_1,\ldots,v_n\}$ $G$ is therefore a chordless $n$-cycle with $n=|G|$. Delete any vertex of $G$; the resulting path on 10 or more vertices is by assumption equivalent to a subgraph of a cyclotomic graph and hence cyclotomic, so a subpath of weight-2 edges is at most 2 edges long. But to exclude graphs of form $X_{4A}$ and $X_{4B}$ the $n$-cycle must contain weight 1 edges only, with one of $v,v'$ (w.l.o.g., $v$) not amongst the $v_i$. Deleting $v$ gives a subgraph with an $n$-cycle that must embed into some cyclotomic graph, so necessarily the cycle is uncharged.

Given the connectivity property, $G$ is therefore \emph{either} of the form:

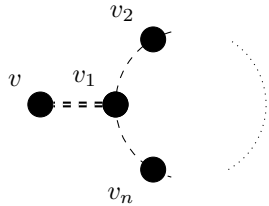
\[
\begin{tikzpicture}[]
\node (a) at (180:1) [fill=black,draw,shape=circle,label=above left:{$v_1$}] {};
\node (c) at (255:1) {};
\draw [dashed] (a) arc(180:240:1);
\node (b1) at (240:1)  {};
\draw [dashed] (b1) arc(-120:-105:1);
\node (b) at (240:1) [fill=black,draw,shape=circle,label=below left:{$v_n$}] {};
\node (d) at (300:1) {};
\draw [dotted] (d) arc(-60:60:1);
\node (f) at (105:1) {};
\draw [dashed] (f) arc (105:180:1);
\node (g) at (120:1) [fill=black,draw,shape=circle,label=above left:{$v_2$}] {};
\node (h) at (-2,0) [fill=black,draw,shape=circle,label=above left:{$v$}] {};
\draw [thick,double,dashed] (a) -- (h);
\node (a2) at (180:1) [fill=black,draw,shape=circle] {};
\end{tikzpicture}
\]
which for $n\ge5$ induces as a proper subgraph on vertices $v,v_1,v_2,v_n,v_{n-1}$ an $\mathcal{L}$-graph equivalent to 
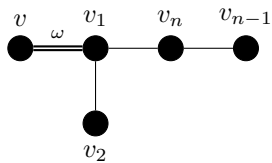
\[
\begin{tikzpicture}
\node (a) at (0,0) [fill=black,draw,shape=circle,label=above:{$v$}] {};
\node (b) at (1,0) [fill=black,draw,shape=circle,label=above:{$v_1$}] {};
\node (c) at (2,0) [fill=black,draw,shape=circle,label=above:{$v_n$}] {};
\node (d) at (3,0) [fill=black,draw,shape=circle,label=above:{$v_{n-1}$}] {};
\node (e) at (1,-1) [fill=black,draw,shape=circle,label=below:{$v_2$}] {};
\draw [double,thick] (a) -- node[above] {{\scriptsize$\omega$}} (b);
\draw (b) -- (c) -- (d);
\draw (b) -- (e);
\end{tikzpicture}
\]
yet no such $\mathcal{L}$-graph is cyclotomic for any $\omega\in\mathcal{L}_2$;

\emph{or}, for some $m\ge1$:
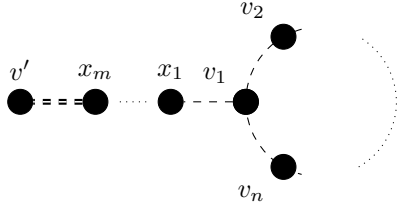
\[
\begin{tikzpicture}[]
\node (a) at (180:1) [fill=black,draw,shape=circle,label=above left:{$v_1$}] {};
\node (c) at (255:1) {};
\draw [dashed] (a) arc(180:240:1);
\node (b1) at (240:1)  {};
\draw [dashed] (b1) arc(-120:-105:1);
\node (b) at (240:1) [fill=black,draw,shape=circle,label=below left:{$v_n$}] {};
\node (d) at (300:1) {};
\draw [dotted] (d) arc(-60:60:1);
\node (f) at (105:1) {};
\draw [dashed] (f) arc (105:180:1);
\node (g) at (120:1) [fill=black,draw,shape=circle,label=above left:{$v_2$}] {};
\node (h) at (-2,0) [fill=black,draw,shape=circle,label=above:{$x_1$}] {};
\draw [dashed] (a) -- (h);
\node (i) at (-3,0) [fill=black,draw,shape=circle,label=above:{$x_m$}] {};
\node (j) at (-4,0) [fill=black,draw,shape=circle,label=above:{$v'$}] {};
\draw [double,thick,dashed] (i) -- (j);
\node (i2) at (-2.8,0) {};
\node (h2) at (-2.2,0) {};
\draw [dotted] (i2) -- (h2);
\node (a2) at (180:1) [fill=black,draw,shape=circle] {};
\end{tikzpicture}
\]
but then the subgraph on vertices $v_1,\ldots,v_n,x_1$ is necessarily a subgraph of some $T_{2k}$, yet this is impossible: if - for a suitable profile - each $v_i\in V_i$ then, as a neighbour of $v_1$, $x_1\in V_2$ or $x_1\in V_n$; yet $x_1$ is not a neighbour of $v_3$ or $v_{n-1}$.  
\end{proof}

\begin{proposition}\label{6327} For $d=-2,-7$ let $G$ be an $\mathcal{L}$-graph with $n\ge11$ vertices, such that every proper connected subgraph of $G$ is equivalent to a subgraph of some $T_{2k}$, $C^{+\pm}_{2k}$, $C^{2+}_{2k}$, $T^4_{2k}$ or ${T^4_{2k}}'$. If $G$ contains an edge label of weight 2, then $G$ is equivalent to a subgraph of some $T_{2k}^4$,$T_{2k}^{4'}$ or $C_{2k}^{2+}$.
\end{proposition}

\begin{proof}
Let $G$ be such a graph: we seek a profile of $G$. Take a chordless path or cycle $P$ with the maximal number of vertices (given a tie, take $P$ to be a path), and let $x$ and $y$ be the endvertices of $P$ if $P$ is a path, or any two adjacent vertices of $P$ if $P$ is a cycle. Note that no vertex of $G$ is adjacent to $x$ but to no other vertex on $P$, else we could either grow $P$ to a longer chordless path, or replace a chordless cycle $P$ by a chordless path of equal length. It follows that $G-\{x\}$ (similarly, $G-\{y\}$) is connected, and since it contains at least 10 vertices it has rank at least 5, so $P$ contains at least 5 vertices. Hence by Proposition \ref{no5cycle} $P$ is necessarily a path, not a cycle.

If there were a vertex not on $P$ adjacent to both $x$ and $y$ but no other vertex on $P$, then $P$ could be extended to a longer chordless cycle, which is impossible. So $G-\{x,y\}$ is connected. It has at least 9 vertices and thus rank $r$ at least 5, so by Lemma \ref{uniqueprofile} it has a uniquely determined profile. As the profiles of $G-\{x\},G-\{y\}$ are also uniquely determined, they can each be obtained by adding $y$ or $x$ to the profile of $G-\{x,y\}$. Since $P$ is not a cycle, $x$ and $y$ are non-adjacent in $G$, and all other possible adjacencies of $x$ in $G$ can be read off from the profile of $G-\{y\}$, and all other possible adjacencies of $y$ in $G$ can be read off from the profile of $G-\{x\}$. Thus we can merge the profiles of $G-\{x\}$ and $G-\{y\}$ to obtain a new sequence of columns $\mathcal{C}$, which we shall show is the profile of $G$. In this merging, columns $2,3,\ldots,r-1$ carry over unchanged, and as $x,y$ are the endpoints of a maximal chordless path they must lie in opposite end columns $1$ and $r$.

Now, no vertex in the column of $x$ is adjacent to one in the column of $y$, else, deleting column 3 of $G-\{x,y\}$ we obtain another proper subgraph of $G$ which thus has a profile that would force all vertices in the column of $x$ to be adjacent to all in the column of $y$. In particular, this would make $x$ a neighbour of $y$ and thus $P$ a cycle. Hence no vertex in column $1$ is adjacent to any in column $r$, and $\mathcal{C}$ is a non-cycling profile of $G$. The local conditions of Proposition \ref{subgraphprop} hold for $G$, since they hold for both $G-\{x\}$ and $G-\{y\}$, so by that result $G$ is equivalent to a subgraph of some $T_{2k},C^{+\pm}_{2k},C^{2+}_{2k},T^4_{2k}$ or $T^{4'}_{2k}$. Since $G$ has at least one edge of weight $2$, it must be contained in one of  $C^{2+}_{2k}$, $T^4_{2k}$ or $T^{4'}_{2k}$.
\end{proof}

\section{Lehmer's Conjecture for $\mathcal{O}_{\Q(\sqrt{d})}$-matrices}

We summarise the results of the preceding Sections to complete the proof of Theorem \ref{lehmer271115}.

\subsection{Proof of Proposition \ref{v10max}}

Proposition \ref{6327} extends the results of \cite{McSm2} to give Proposition \ref{newprop8}. Thus if $G$ is a non-supersporadic $\mathcal{L}'$-graph with at least eleven vertices, then $G$ is contained in a cyclotomic graph and therefore cyclotomic itself. A minimal noncyclotomic $\mathcal{L}'$-graph therefore either has at most ten vertices, or is supersporadic. But by the results of section \ref{supersporadicsection}, the latter case also forces $G$ to have at most ten vertices.

\subsection{Lehmer's Conjecture for $\mathcal{L}'$-graphs}
We conclude the following:

\begin{theorem}\label{Ldashlehmer} For $d=-2,-7$, if $G$ is a noncyclotomic $\mathcal{L}'$-graph, then $M(G)\ge\lambda_0$.
\end{theorem}

\begin{proof} If $G$ is noncyclotomic then it contains a minimal noncyclotomic subgraph $G'$. By Proposition \ref{v10max}, $G'$ has at most ten vertices. If $G'$ has a weight 2 edge, then by Proposition \ref{nosmallsmall} $M(G)>1.3>\lambda_0$. But if $G'$ is an $\mathcal{L}'$-graph without a weight 2 edge, then it is an $(\mathcal{L}_1\cup\{0\})$-graph; for $d=-2,-7$ we have $\mathcal{L}_1=\{\pm1\}$, and thus $G'$ is a charged signed graph. But then by \cite{McSm2} $M(G')\ge\lambda_0$. In either case, we therefore have $M(G)\ge M(G')\ge\lambda_0$.
\end{proof}

\subsection{Proof of Theorem \ref{lehmer271115} for $d=-2$}
Let $A$ be a noncyclotomic $\mathcal{O}_{\Q(\sqrt{-2})}$-matrix such that $\displaystyle\max_{i}|a_{ii}|=m$ and $\displaystyle\max_{i\neq j} \mbox{Norm}(a_{ij})=n$. If $m\ge 3$ then $M(A)\ge 2.618$ by Proposition \ref{diag2max}. If $m\le 2$ but $n\ge 5$ then $M(A)\ge 2.36$ by Proposition \ref{offdiag4max}. If $n\le 4$ but $m=2$ then $M(A)\ge 1.722$ by Proposition \ref{diag2max2}. Otherwise, $A$ is an $\mathcal{L}$-matrix with $\mathcal{L}$-graph representative $G$. If $n=4$ $M(A)=M(G)\ge 2.08$ by Proposition \ref{w4edges}; whilst if $n=3$ then $M(A)=M(G)\ge 1.56$ by Proposition \ref{w3edges}. But if not, then $G$ is an $\mathcal{L}'$-graph and thus $M(A)\ge\lambda_0$ by Theorem \ref{Ldashlehmer}.

\subsection{Proof of Theorem \ref{lehmer271115} for $d=-7$}
Let $A$ be a noncyclotomic $\mathcal{O}_{\Q(\sqrt{-7})}$-matrix such that $\displaystyle\max_{i}|a_{ii}|=m$ and $\displaystyle\max_{i\neq j} \mbox{Norm}(a_{ij})=n$. If $m\ge 3$ then $M(A)\ge 2.618$ by Proposition \ref{diag2max}. If $m\le 2$ but $n\ge 5$ then $M(A)\ge 2.36$ by Proposition \ref{offdiag4max}. If $n\le 4$ but $m=2$ then $M(A)\ge 1.722$ by Proposition \ref{diag2max2}. Otherwise, $A$ is an $\mathcal{L}$-matrix with $\mathcal{L}$-graph representative $G$. If $n=4$ $M(A)=M(G)\ge 2.08$ by Proposition \ref{w4edges}; otherwise (since $\mathcal{L}_3=\emptyset$) $G$ is an $\mathcal{L}'$-graph and thus $M(A)\ge\lambda_0$ by Theorem \ref{Ldashlehmer}.

\subsection{Proof of Theorem \ref{lehmer271115} for $d=-11$}
Let $A$ be a noncyclotomic $\mathcal{O}_{\Q(\sqrt{-11})}$-matrix such that $\displaystyle\max_{i}|a_{ii}|=m$ and $\displaystyle\max_{i\neq j} \mbox{Norm}(a_{ij})=n$. If $m\ge 3$ then $M(A)\ge 2.618$ by Proposition \ref{diag2max}. If $m\le 2$ but $n\ge 5$ then $M(A)\ge 2.36$ by Proposition \ref{offdiag4max}. If $n\le 4$ but $m=2$ then $M(A)\ge 1.722$ by Proposition \ref{diag2max2}. Otherwise, $A$ is an $\mathcal{L}$-matrix with $\mathcal{L}$-graph representative $G$. If $n=4$ $M(A)=M(G)\ge 2.08$ by Proposition \ref{w4edges}; whilst if $n=3$ then $M(A)=M(G)\ge 1.56$ by Proposition \ref{w3edges}. Otherwise (since $\mathcal{L}_2=\emptyset$) $G$ is a charged signed graph, so $M(A)=M(G)\ge\lambda_0$ by \cite{McSm2}.

\subsection{Proof of Theorem \ref{lehmer271115} for $d=-15$}
Let $A$ be a noncyclotomic $\mathcal{O}_{\Q(\sqrt{-15})}$-matrix such that $\displaystyle\max_{i}|a_{ii}|=m$ and $\displaystyle\max_{i\neq j} \mbox{Norm}(a_{ij})=n$. If $m\ge 3$ then $M(A)\ge 2.618$ by Proposition \ref{diag2max}. If $m\le 2$ but $n\ge 5$ then $M(A)\ge 2.36$ by Proposition \ref{offdiag4max}. If $n\le 4$ but $m=2$ then $M(A)\ge 1.722$ by Proposition \ref{diag2max2}. Otherwise, $A$ is an $\mathcal{L}$-matrix with $\mathcal{L}$-graph representative $G$. If $n=4$ $M(A)=M(G)\ge 2.08$ by Proposition \ref{w4edges}. Otherwise (since $\mathcal{L}_3=\mathcal{L}_2=\emptyset$) $G$ is a charged signed graph, so $M(A)=M(G)\ge\lambda_0$ by \cite{McSm2}.

\subsection{Proof of Theorem \ref{lehmer271115} for other $d$}
Let $A$ be a noncyclotomic $\mathcal{O}_{\Q(\sqrt{d})}$-matrix for squarefree $d\le -17$ or $d\in\{-5,-6,-10,-13,-14\}$ such that $\displaystyle\max_{i}|a_{ii}|=m$ and $\displaystyle\max_{i\neq j} \mbox{Norm}(a_{ij})=n$. If $m\ge 3$ then $M(A)\ge 2.618$ by Proposition \ref{diag2max}. If $m\le 2$ but $n\ge 5$ then $M(A)\ge 2.36$ by Proposition \ref{offdiag4max}. Otherwise (since $\mathcal{L}=\{0,\pm 1,\pm 2\}$) $A$ is an integer symmetric matrix so $M(A)\ge\lambda_0$ by \cite{McSm2}.

\section{Acknowledgements}
This work has made use of the resources provided by the Edinburgh Computer and Data Facility (ECDF). (\url{http://www.ecdf.ed.ac.uk}). The ECDF is partially supported by the eDIKT initiative (\url{http://www.edikt.org.uk}).

All tests for cyclotomicity / minimal noncyclotomicity were performed using the SAGE computer algebra system (\cite{sage}); larger calculations, including the search described in Section 4.4, were performed in parallel using the ECDF. For source code and other implementation details, see \cite{Ta1}.

\begin{spacing}{0.4}

\end{spacing}

\end{document}